\documentclass{elsart}
\journal{J. Number Theory}
\usepackage{amssymb}
\usepackage{amsmath}

\newtheorem{theorem}{Theorem}[section]
\newtheorem{lemma}{Lemma}[section]
\newenvironment{proof}%
{\rm \trivlist \item[\hskip \labelsep{\bf Proof. }]}%
{\hspace*{\fill}$\Box$\endtrivlist}
\numberwithin{equation}{section}

\DeclareRobustCommand{\qbinom}{\genfrac[]{0pt}{}}

\begin{document}
\begin{frontmatter}
\title{Irrationality of $\zeta_q(1)$ and $\zeta_q(2)$ \thanksref{label1}}
\thanks[label1]{This work was supported by INTAS Research Network NeCCA (03-51-6637),
FWO project G.0455.04 and OT/04/21 of K.U.Leuven.}
\author{Kelly Postelmans and} 
\author{Walter Van Assche}
\ead{Walter.VanAssche@wis.kuleuven.be}

\address{Katholieke Universiteit Leuven, Department of Mathematics,\\
Celestijnenlaan 200B, B-3001 Leuven, Belgium}

\begin{abstract}
In this paper we show how one can obtain simultaneous rational
approximants for $\zeta_q(1)$ and $\zeta_q(2)$ with a common
denominator  by means of Hermite-Pad\'e approximation using
multiple little $q$-Jacobi polynomials and we show that properties
of these rational approximants prove that $1$, $\zeta_q(1)$,
$\zeta_q(2)$ are linearly independent over $\mathbb{Q}$. In
particular this implies that $\zeta_q(1)$ and $\zeta_q(2)$ are
irrational. Furthermore we give an upper bound for the measure of
irrationality.
\end{abstract}

\begin{keyword}
$q$-zeta function, irrationality, simultaneous rational approximation
\end{keyword}
\end{frontmatter}

\section{Introduction}
There are a lot of things known about the irrationality of the
Riemann zeta-function
\[     \zeta(s)=\sum_{n=1}^{\infty} \frac{1}{n^s}  \]
at integer values $s \in \{2,3,4,\ldots\}$.  It is known that
$\zeta(2n)=(-1)^{n-1}2^{2n-1}B_{2n} \pi^{2n}/(2n)!$, where
$B_{2n}$ are Bernoulli numbers (which are rational), so it follows
that $\zeta(2n)$ is an irrational number for $n \in
\mathbb{N}_0=\{1,2,3,\ldots\}$. It is also known that $\zeta(3)$
is irrational (Ap\'ery \cite{Ap}) and that $\zeta(2n+1)$ is irrational for
an infinite number of $n \in \mathbb{N}_0$ \cite{BallRivoal}
\cite{Fischler} \cite{Zudilin}. Furthermore at least one of the
numbers $\zeta(5), \zeta(7), \zeta(9), \zeta(11)$ is irrational
\cite{Rivoal} \cite{Zudilin3}.

A possible way for a $q$-extension of the  Riemann zeta-function is
\begin{equation}\label{qzetas}
    \zeta_q(s)=\sum_{n=1}^\infty \frac{n^{s-1}q^n}{1-q^n}, \qquad
    s=1,2, \ldots,
\end{equation}
with $q \in \mathbb{C}$ and $|q|<1$. Then the limit relations
\begin{equation}
    \lim_{\substack{q\to1\\|q|<1}}(1-q)^s\zeta_q(s)=(s-1)! \,
    \zeta(s), \qquad s=2,3, \ldots,
\end{equation}
hold \cite{Zudilin2}. Earlier it was shown that $\zeta_q(1)$ is
irrational whenever $q=1/p$, with $p$ an integer greater than 1
(see, e.g., \cite{VanAssche} and the references there). Results of
Nesterenko \cite{Nesterenko} show that $\zeta_q(2)$ is
transcendental for every algebraic number $q$ with $0< |q|<1$.
Zudilin gave an upper bound for the measure of
irrationality of $\zeta_q(2)$ \cite{Zudilin1} with $1/q \in
\{2,3,4,\ldots\}$. Furthermore Krattenthaler, Rivoal and Zudilin \cite{KRZ}
proved that there are infinitely many $\zeta_q(2n)$ (and infinitely many
$\zeta_q(2n+1)$) which are irrational, and that at least one of the numbers
$\zeta_q(3),\zeta_q(5),\zeta_q(7),\zeta_q(9),\zeta_q(11)$ is irrational
whenever $1/q \in \mathbb{Z}$ and $q \neq \pm 1$.
From now we only use values of $q$ for which
$p=1/q \in \mathbb{N}\setminus\{0,1\} = \{2,3,4,\ldots\}$.

In this paper we will show how one can obtain good
\textit{simultaneous} rational approximations for $\zeta_q(1)$ and
$\zeta_q(2)$ with a common denominator, which are related to
multiple little $q$-Jacobi polynomials and Hermite-Pad\'e
approximation techniques. The method that we use is an extension
of the method that Van Assche \cite{VanAssche}
 used to prove the
irrationality of $\zeta_q(1)$ and is an application of
Hermite-Pad\'e approximation of a system of Markov functions
\cite[Chapter 4]{NikiSor} \cite{VanAssche2}.

In Section 2 we will describe the construction of the simultaneous
rational approximants using Hermite-Pad\'e approximation to two
Markov functions which are chosen appropriately. The solution of
this Hermite-Pad\'e approximation problem depends on multiple
little $q$-Jacobi polynomials \cite{Walter-Kelly}, and in Section
3 we give the relevant formulas for the special case where all
parameters are equal to zero. In Section 4 we will show that these
rational approximations give a good rational approximation to
$\zeta_q(1)$ and in particular we will prove the following result:
\begin{theorem} \label{thm:qzeta1}
Suppose $p>1$ is an integer and $q=1/p$. Let
$\alpha_n$ and $\beta_n$ be given by (\ref{beta})--(\ref{alpha})
for all $n \in \mathbb{N}$. Then $\alpha_n, \beta_n \in
\mathbb{Z}$ and $\beta_n \zeta_q(2) - \alpha_n \neq 0$.
Furthermore
\[   \lim_{n \to \infty} \left|\beta_n \zeta_q(1) - \alpha_n\right|^{1/n^2}
    \leq p^{-\frac{3(\pi^2-4)}{\pi^2}}<1. \]
\end{theorem}

In Section 5 we show that we also get a good rational
approximation to $\zeta_q(2)$:
\begin{theorem} \label{thm:qzeta2}
Suppose $p>1$ is an integer and $q=1/p$. Let 
$a_n$ and $b_n$ be given by (\ref{bnm})--(\ref{anm}) for all $n \in
\mathbb{N}$. Then $a_n, b_n \in \mathbb{Z}$ and $b_n \zeta_q(2) -
a_n \neq 0$. Furthermore
\[  \lim_{n \to \infty} \left|b_n \zeta_q(2) - a_n\right|^{1/n^2}
    \leq p^{-\frac{3(\pi^2-8)}{\pi^2}}<1.
\]
\end{theorem}
These rational approximations are good enough to conclude that
$\zeta_q(1)$ and $\zeta_q(2)$ are irrational, by using the
following elementary lemma:

\begin{lemma}\label{lem:irrational}
Let $x$ be a real number and suppose there exist integers $a_n,
b_n$ $(n \in \mathbb{N})$ such that
\begin{enumerate}
    \item $b_nx-a_n \neq 0$ for every $n \in \mathbb{N}$,
    \item $\lim_{n\to \infty} |b_n x - a_n|=0$.
\end{enumerate}
Then $x$ is irrational.
\end{lemma}
This lemma expresses the fact that a rational number is
approximable to order 1 by rational numbers and to no higher order
\cite[Theorem 186]{Hardy-Wright}. The measure of irrationality
$r(x)$ can be defined as
\[     r(x)=\inf\{r \in \mathbb{R}: |x-a/b|<1/b^r \textrm{ has at most
    finitely many integer solutions } (a,b) \}. \]
It is known that if $|x-a_n/b_n|=\mathcal{O}(1/b_n^{1+s})$ with $0<s<1$ and
$b_n<b_{n+1}<b_n^{1+o(1)}$, then the measure of irrationality
$r(x)$ satisfies $2 \leq r(x) \leq 1+1/s$ (see, e.g.,
\cite[exercise 3 on p.~376]{Borwein-Borwein} for the upper bound;
the lower bound follows since every irrational number is
approximable to order 2 \cite[Theorem 187]{Hardy-Wright}).

Our rational approximations to $\zeta_q(1)$ give the upper bound
$5.0444\ldots$ for the measure of irrationality of $\zeta_q(1)$,
which is not as good as the upper bound of $2.5082\ldots$ in
\cite{VanAssche} or the upper bound of $2.4234\ldots$ in
\cite{Zudilin2}. For $\zeta_q(2)$ we obtain the upper bound
$15.8369\ldots$ for the measure of irrationality, which is not as
good as the upper bound $4.07869374\ldots$ that Zudilin obtained
in \cite{Zudilin1}. Our rational approximations, however, use the
same denominator and hence we have constructed simultaneous
rational approximants. The price we pay for this is that the order
of approximation for each individual number is not as good as
possible. But we gain some very important information because our
simultaneous rational approximants are good enough to prove our
main result:

\begin{theorem} \label{thm:Linind}
The numbers $1$, $\zeta_q(1)$ and $\zeta_q(2)$ are linearly
independent over $\mathbb{Q}$.
\end{theorem}

This result is  stronger than the statement that $\zeta_q(1)$ and
$\zeta_q(2)$ are irrational. We prove this theorem by means of the
following lemma, which extends Lemma \ref{lem:irrational} (see
\cite[Lemma 2.1]{Hata} for a similar lemma but with a different
first condition).

\begin{lemma}\label{lem:Linind}
Let $x, y $ be real numbers. Suppose that  for all $(a, b, c) \in
\mathbb{Z}^3\setminus(0,0,0)$ and for infinity many positive
integers $n \in \Lambda \subset \mathbb{N}$, there exist integers
$p_n$, $q_n$ and $r_n$ such that the following three conditions
are satisfied:
\begin{enumerate}
 \item $a p_n+b q_n+c r_n \neq 0$ for every $n \in \Lambda$,
 \item $|p_n x - q_n| \to 0$ as $n \to \infty$,
 \item $|p_n y - r_n| \to 0$ as $n \to \infty$.
\end{enumerate}
Then $1, x, y$ are linearly independent over $\mathbb{Q}$.
\end{lemma}

\begin{proof}
Suppose that $1,x$ and $y$ are linearly dependent over
$\mathbb{Q}$. Then there exist integers $a_1$, $a_2$, $b_1$,
$b_2$, $c_1$, $c_2$ such that
\[     \frac{a_1}{a_2} + \frac{b_1}{b_2}x + \frac{c_1}{c_2}y=0. \]
When we multiply this by $a_2b_2c_2$ we get $a_1b_2c_2 +
a_2b_1c_2x + a_2b_2c_1y=0$, hence there exist integers $a,b,c$
such that $a+bx+cy=0$. Multiplying by $p_n$ and adding the terms
$-bq_n-cr_n$ on both sides, we obtain
\[     b(p_nx-q_n)+c(p_ny-r_n)=-(ap_n+bq_n+cr_n). \]
On the right hand side we have an integer different from zero,
hence the right hand side is in absolute value at least 1 for
every $n \in \Lambda$. But the expression on the left hand side
tends to $0$ as $n$ tends to infinity, hence we have a
contradiction.
\end{proof}

If we take $x=\zeta_q(1)$ and $y= \zeta_q(2)$, then Condition (2)
follows from Theorem \ref{thm:qzeta1} and Condition (3) follows from
Theorem \ref{thm:qzeta2}. In Section 6 we will show that Condition
(1) of Lemma \ref{lem:Linind} holds. This requires some results
about cyclotomic polynomials and cyclotomic numbers from number
theory.

\section{Multiple little $q$-Jacobi polynomials}
\subsection{Multiple little $q$-Jacobi polynomials}
 Little $q$-Jacobi polynomials are orthogonal
polynomials on the exponential lattice $\{q^k,\ k=0,1,2,\ldots\}$,
where $0 < q < 1$. In order to express the orthogonality
relations, we will use the $q$-integral
\begin{equation} \label{eq:qint}
   \int_0^1 f(x)\, d_qx =  \sum_{k=0}^\infty q^k f(q^k),
\end{equation}
(see, e.g., \cite[\S 10.1]{AAR}, \cite[\S 1.11]{GR}) where $f$ is
a function on $[0,1]$ which is continuous at $0$. The
orthogonality is given by
\begin{equation} \label{eq:ortho}
  \int_0^1 p_n(x;\alpha,\beta|q) x^k w(x;\alpha,\beta|q)\, d_qx = 0, \qquad k=0,1,\ldots,n-1,
\end{equation}
where
\begin{equation} \label{eq:w}
   w(x;a,b|q) = \frac{(qx;q)_\infty}{(q^{\beta+1}x;q)_\infty} x^\alpha .
\end{equation}
We have used the notation
\[  (a;q)_n = \prod_{k=0}^{n-1} (1-aq^k), \qquad (a;q)_\infty = \prod_{k=0}^\infty (1-aq^k). \]
In order that the $q$-integral of $w$ is finite, we need to impose
the restrictions $\alpha,\beta > -1$. The orthogonality conditions
(\ref{eq:ortho}) determine the polynomials $p_n(x;\alpha,\beta|q)$
up to a multiplicative factor. When we define
$p_n(x;\alpha,\beta|q)$ as monic polynomials, they are uniquely
determined by the orthogonality conditions.

The Rodrigues formula for $p_n(x;\alpha,\beta|q)$ is given by
\begin{equation}\label{rodrigues1}
  p_n(x;\alpha,\beta|q)=
    \frac{(q^{\alpha+n}-q^{\alpha+n-1})^n}
    {(q^{\alpha+\beta+n+1};q)_n w(x,\alpha, \beta|q)}
    D_{q^{-1}}^n w(x, \alpha+n, \beta + n|q)
\end{equation}
and an explicit formula by
\begin{eqnarray}
  p_n(x;\alpha,\beta|q) &=& (-1)^n q^{n(n-1)/2}
    \frac{(q^{\alpha+1};q)_n}{(q^{\alpha+\beta+n+1};q)_n}
   {}_2\phi_1 \left( \left. \begin{array}{c}
   q^{-n},q^{n+\alpha+\beta+1} \\ q^{\alpha+1}
   \end{array} \right| q;qx \right) \nonumber \\
  & = & (-1)^n q^{n(n-1)/2}
    \frac{(q^{\alpha+1};q)_n}{(q^{\alpha+\beta+n+1};q)_n}
  \sum_{k=0}^n \frac{(q^{-n};q)_k (q^{n+\alpha+\beta+1};q)_k}{(q^{\alpha+1};q)_k} \frac{q^kx^k}{(q;q)_k} .  \label{eq:expl}
\end{eqnarray}

We have introduced multiple little $q$-Jacobi polynomials in \cite{Walter-Kelly}. Suppose that 
$\beta>-1$ and that $\alpha_1, \alpha_2, \ldots, \alpha_r$ are such that
each $\alpha_i>-1$ and $\alpha_i-\alpha_j \not\in \mathbb{Z}$
whenever $i\neq j$. Then the  multiple little $q$-Jacobi
polynomial $p_{\vec{n}}(x;\vec{\alpha},\beta|q)$ for the
multi-index $\vec{n}=(n_1, n_2, \ldots, n_r)\in \mathbb{N}^r$ is
the monic polynomial of degree $|\vec{n}|=n_1+n_2+ \ldots+n_r$
that satisfies the orthogonality conditions
\begin{equation} \label{eq:mlqJ}
  \int_0^1 p_{\vec{n}}(x;\vec{\alpha},\beta|q) x^k w(x;\alpha_j,\beta|q)\, d_qx  =  0, \qquad k=0,1,\ldots,n_j-1,\ j=1,2,\ldots,r,
\end{equation}
These are the multiple little $q$-Jacobi polynomials of the first
kind. We only consider the case $r=2$, so we take $\vec{n}=(n,m)$,
with $m \leq n$. From \cite{Walter-Kelly} we know that the
Rodrigues formula for $p_{n,m}(x;(\alpha_1, \alpha_2),\beta|q)$ is
given by
\begin{multline}\label{rodrigues}
    p_{n,m}(x;(\alpha_1, \alpha_2),\beta|q)=
    \frac{(q^{\alpha_1+n+m}-q^{\alpha_1+n+m-1})^n (q^{\alpha_2+m}-q^{\alpha_2+m-1})^m}
    {(q^{\beta+\alpha_1+n+m+1};q)_n (q^{\beta+\alpha_2+n+m+1};q)_m w(x, \alpha_1,
    \beta|q)}\\
     \times  D_{q^{-1}}^n[x^{\alpha_1-\alpha_2+n}
    D_{q^{-1}}^m w(x, \alpha_2+m, \beta+n+m|q)]
\end{multline}
and that an explicit expression for $p_{n,m}(x;(\alpha_1,
\alpha_2),\beta|q)$ is given by
\begin{multline}\label{uitdr mlqJ}
    p_{n,m}(x;(\alpha_1, \alpha_2),\beta|q)= (-1)^{n+m}
    \frac{ q^{nm} (q^{\alpha_1 +1};q)_n (q^{\alpha_2+1};q)_m q^{n(n-1)/2} q^{m(m-1)/2}}
    {(q^{\beta+\alpha_1+n+m+1};q)_n
    (q^{\beta+\alpha_2+n+m+1};q)_m}\\
    \times
    \sum_{k=0}^m \sum_{j=0}^n
    \frac{(q^{-n};q)_j(q^{-m};q)_k (q^{\alpha_1+n+1};q)_k(q^{\alpha_1+\beta+n+1};q)_{j+k}
    (q^{\alpha_2+\beta+n+m+1};q)_k}
    {(q^{\alpha_1+1};q)_{j+k}
    (q^{\alpha_1+\beta+n+1};q)_k
    (q^{\alpha_2+1};q)_k}\\
    \times
    \frac{q^{k+j} x^{k+j}}{q^{kn}(q;q)_k(q;q)_j}.
\end{multline}

\subsection{The case $\alpha_1=\alpha_2$}
When we take $\alpha_1=\alpha_2= \alpha$, then (\ref{eq:mlqJ}) gives 
only $n$ conditions. However, when we look at expression
(\ref{uitdr mlqJ}) we see that $p_{n,m}(x;(\alpha_1,
\alpha_2),\beta|q)$ still has degree $n+m$. We use the notation
$p_{n, m}^{(\alpha, \alpha,
        \beta)}(x)$
 for $p_{n,m}(x;(\alpha_1, \alpha_2),\beta|q)$ with $\alpha_1=\alpha_2= \alpha$.
We will now show that $p_{n, m}^{(\alpha, \alpha,
        \beta)}(x)$, defined by the Rodrigues formula
(\ref{rodrigues}), again satisfies $n+m$ orthogonality conditions, namely
\begin{eqnarray}           
     \int_0^1 p_{n, m}^{(\alpha, \alpha,
        \beta)}(x) w(x, \alpha,
        \beta|q) x^\ell\, d_qx & = & 0,
     \qquad \ell=0,1, \ldots, n-1,  \label{vwpnab1}  \\  
    \int_0^1 p_{n, m}^{(\alpha, \alpha, 
        \beta)}(x) w(x, \alpha,
        \beta|q) x^\ell \log_qx \, d_qx & = & 0, 
     \qquad \ell=0,1, \ldots, m-1. \label{vwpnab2}
 \end{eqnarray}
To prove these orthogonality conditions we will need the following
lemma.
\begin{lemma}[Summation by parts]
When $g(q^{-1})=0$ we have that
\begin{equation}\label{partsomp_q}
     \sum_{k=0}^\infty q^k f(q^k) \left. D_{q^{-1}} g(x) \right|_{q^k}
    = -q  \sum_{k=0}^\infty q^k g(q^k) \left. D_{q} f(x) \right|_{q^k} .
\end{equation}
\end{lemma}
We start with the first $n$ orthogonality conditions. When we use
the Rodrigues formula (\ref{rodrigues}), we obtain
\[
\int_0^1 p_{n, m}^{(\alpha, \alpha, \beta)}(x) w(x, \alpha,
\beta|q) x^\ell\, d_qx =   C_{n,m} \int_0^1
      D_{q^{-1}}^n  \left[x^{n} D_{q^{-1}}^m
     w(x, \alpha+m, \beta+n+m|q)\right]
  x^\ell\, d_qx,
\]
with 
\[ C_{n,m}=\frac{(q^{\alpha+n+m}-q^{\alpha+n+m-1})^n
(q^{\alpha+m}-q^{\alpha+m-1})^m}{(q^{\beta+\alpha+n+m+1};q)_n
    (q^{\beta+\alpha+n+m+1};q)_m} .  \]
Now we can rewrite the right hand side as an infinite sum using
(\ref{eq:qint}) and then apply summation by parts
(\ref{partsomp_q}) $n$ times. This gives
\begin{multline*}
\int_0^1 p_{n, m}^{(\alpha, \alpha, \beta)}(x) w(x, \alpha,
\beta|q) x^\ell\, d_qx \\ 
   =  (-1)^n q^n C_{n,m} \sum_{k=0}^\infty q^k
        \left.  \rule{0pt}{15pt} \left[x^{n} D_{q^{-1}}^m
     w(x, \alpha+m, \beta+n+m|q) \right]
     \right|_{x=q^k} \left. \rule{0pt}{15pt}
  D_q^n x^\ell \right|_{x=q^k}.
\end{multline*}
Since $D_{q}^n x^\ell=0$ for $\ell < n$, we see that (\ref{vwpnab1}) holds.

Now we prove (\ref{vwpnab2}). Again using
the Rodrigues formula (\ref{rodrigues}) and applying summation by
parts (\ref{partsomp_q}) $n$ times 
\begin{multline*}
    \int_0^1 p_{n, m}^{(\alpha, \alpha, \beta)}(x) 
    w(x, \alpha, \beta|q) x^\ell \log_qx \, d_qx \\ 
   = (-1)^n  q^n C_{n,m}   \int_0^1 D_{q}^n  \left( x^\ell \log_qx   \right) 
   \left[ x^n  D_{q^{-1}}^m  w(x, \alpha+m, \beta+n+m|q) \right]  \, d_qx.
\end{multline*}
To calculate $D_{q}^n  \left( x^\ell \log_qx   \right)$ we can use
the $q$-analogue of Leibniz' rule
\begin{equation}\label{Leibniz}
    D_q^n[f(x)g(x)]=\sum_{k=0}^n
    \qbinom{n}{k}_q \left( D_q^{n-k} f
    \right) (q^k x) \left(D_q^k g\right)(x), \quad n=0,1,2,\ldots.
\end{equation}
This gives
\begin{multline*}
    \int_0^1 p_{n, m}^{(\alpha, \alpha, \beta)}(x) 
    w(x, \alpha,\beta|q) x^\ell \log_qx \, d_qx  \\ 
    = C   \int_0^1 \sum_{k=0}^n
    \qbinom{n}{k}_q \left. \rule{0pt}{15pt}\left( D_{q}^{n-k}  x^\ell  \right)
   \right|_{q^{k}x}  \left. \rule{0pt}{15pt} \left(D_{q}^k \log_q
    x\right)\right|_x
    \left[ x^{n}  D_{q^{-1}}^m  w(x, \alpha+m, \beta+n+m|q) \right]  \, d_qx,
\end{multline*}
with $C$ a constant. Obviously we have for $k <n-\ell$ that $D_{q}^{n-k} x^\ell=0$ 
and for  $k \geq n-\ell$ that $ D_{q}^{n-k} x^\ell$ is a polynomial of degree 
$\ell-n+k$, say $\pi_{k,\ell-n+k}$. It is also
easy to see that $D_{q}^k \log_q x$ is a constant times
$1/x^k$ for $k \geq 1$. So we get
\begin{multline*}
    \int_0^1 p_{n, m}^{(\alpha, \alpha,\beta)}(x) 
    w(x, \alpha, \beta|q) x^\ell \log_qx \, d_qx  \\  
    =  \int_0^1 \sum_{k=n-\ell}^n c_k
    \qbinom{n}{k}_q
    \pi_{k,\ell-n+k}(q^kx)\frac{1}{x^k}
    x^{n}   D_{q^{-1}}^m  w(x, \alpha+m, \beta+n+m|q) \, d_qx,
\end{multline*}
where the $c_k$ are constants. Now $\pi_{k,\ell-n+k}(q^{k}x) x^{n}/x^k$ is a
polynomial of degree $\ell$, say $\rho_{k,\ell}(x)$. 
Using the Rodrigues formula (\ref{rodrigues1}) for little $q$-Jacobi polynomials, we get
\begin{multline*}
    \int_0^1  p_{n, m}^{(\alpha, \alpha, \beta)}(x) 
    w(x, \alpha, \beta|q) x^\ell \log_qx \, d_qx \\  
   = (-1)^n  \sum_{k=n-\ell}^n c'_k  \qbinom{n}{k}_q
    \int_0^1 \rho_{k,\ell}(x)
     w(x, \alpha, \beta+n|q) p_m^{(\alpha, \beta+n)}(x) \, d_qx,
\end{multline*}
where the $c'_k$ are constants. The orthogonality relations of
the little $q$-Jacobi polynomial imply that the right hand side
is $0$ for $\ell<m$, so that (\ref{vwpnab2}) holds.

\subsection{The case $\alpha_1=\alpha_2=\beta=0$}
From now, we will work with the multiple little $q$-Jacobi
polynomial $p_{n,m}(x;(\alpha_1, \alpha_2),\beta|q)$ where
$\alpha_1=\alpha_2=\beta=0$ and with
    another normalisation, namely
\begin{equation}\label{norm2}
    p_{n,m}(0;(\alpha_1,\alpha_2),\beta|q)=1.
\end{equation}
We will denote this polynomial by $p_{n, m}(x)$. Since
$w(x,0,0|q)=1$, it follows from (\ref{vwpnab1})--(\ref{vwpnab2}) that $p_{n,m}$ can
be defined by
\begin{eqnarray}  
     \int_0^1 p_{n, m}(x) x^\ell\, d_qx & = & 0, 
     \qquad \ell=0,1, \ldots, n-1,   \label{vwpn1} \\
    \int_0^1 p_{n, m}(x) x^\ell \log_qx \, d_qx & = & 0,
    \qquad \ell=0,1, \ldots, m-1,  \label{vwpn2}
 \end{eqnarray}
together with the normalisation (\ref{norm2}).

 From (\ref{uitdr mlqJ}) we  see that
\begin{equation}\label{kappainverse}
    p_{n,m}(0;(\alpha_1,\alpha_2),\beta|q)   = (-1)^{n+m}
    \frac{ q^{nm} (q^{\alpha_1 +1};q)_n (q^{\alpha_2+1};q)_m q^{n(n-1)/2} q^{m(m-1)/2}}
    {(q^{\beta+\alpha_1+n+m+1};q)_n
    (q^{\beta+\alpha_2+n+m+1};q)_m}.
\end{equation}
Multiplying (\ref{rodrigues}) and (\ref{uitdr mlqJ}) by
$[p_{n,m}(0;(\alpha_1, \alpha_2),\beta|q)]^{-1}$ and setting
$\alpha_1=\alpha_2=\beta=0$ gives
\begin{equation}\label{mlqJa}
 p_{n, m}(x)=
    \sum_{k=0}^m \sum_{j=0}^n \frac{(qx)^{k+j}}{(q;q)_k(q;q)_j}
    \frac{1}{q^{kn}}
     \frac{(q^{n+1};q)_{j+k}}{(q;q)_{j+k}}
    \frac{(q^{n+m+1};q)_k}{(q;q)_k} (q^{-n};q)_j(q^{-m};q)_k
\end{equation}
and the Rodrigues formula becomes
\begin{equation}\label{rodriguesa}
 p_{n, m}(x)=(-1)^{n+m} q^{n(n-1)/2} q^{m(m-1)/2} \frac{(q-1)^{n+m}}
    {(q;q)_n (q;q)_m} D_{q^{-1}}^n[x^{n}
    D_{q^{-1}}^m \left((qx;q)_{n+m} x^m \right)].
\end{equation}
Equation (\ref{mlqJa}) expresses the polynomial $p_{n,m}$
in the basis $\{1, x, x^2, \ldots, x^{n+m}\}$. Sometimes it is
more convenient to use the basis $\{(qx;q)_\ell, 0 \leq \ell \leq n+m
\}$. The Rodrigues formula
(\ref{rodriguesa}) allows us to obtain an expression for
$p_{n, m}(x)$ in this basis. Recall the $q$-analog of Newton's
binomial formula
\begin{equation}\label{artw12}
    (x;q)_n=  \sum_{k=0}^n
    \qbinom{n}{k}_q q^{k(k-1)/2}(-x)^k,
\end{equation}
and its dual
\begin{equation}\label{artw13}
    x^n= \sum_{k=0}^n
    \qbinom{n}{k}_q  (-1)^k q^{-nk +
    k(k+1)/2}(x;q)_k.
\end{equation}
This is a special case of the $q$-binomial series \cite[\S 10.2]{AAR} \cite[\S 1.3]{GR}
\begin{equation}\label{artw14}
    \sum_{n=0}^\infty\frac{(a;q)_n}{(q;q)_n}x^n=
    \frac{(ax;q)_\infty}{(x;q)_\infty}, \quad \mbox{$|q|<1$, $|x|<1$.}
\end{equation}
Using (\ref{artw13}) with argument $q^{n+m+1}x$ and exponent $m$
gives 
\[
    x^m= p^{(n+m+1)m}\sum_{k=0}^m
    \qbinom{m}{k}_q  (-1)^k q^{-mk +
    k(k+1)/2}(q^{n+m+1}x;q)_k.
\]
Using this in the Rodrigues formula (\ref{rodriguesa}), we
find
\begin{multline*}
    p_{n, m}(x)=
    (-1)^{n+m} q^{n(n-1)/2} q^{m(m-1)/2} \frac{(q-1)^{n+m}}
    {(q;q)_n (q;q)_m}
    p^{(n+m+1)m}\\
    \sum_{k=0}^m
    \qbinom{m}{k}_q  (-1)^k q^{-mk +
    k(k+1)/2}
     D_{q^{-1}}^n[x^{n}
    D_{q^{-1}}^m (qx;q)_{n+m+k} ],
\end{multline*}
because $(qx;q)_{n+m}\, (q^{n+m+1}x;q)_k=(qx;q)_{n+m+k}$. One easily
finds that
\begin{equation}\label{artw15}
    D_p^k(qx;q)_n= \frac{(q;q)_n}{(q;q)_{n-k}(1-p)^k}(qx;q)_{n-k},
\end{equation}
so we get
\begin{multline*}
    p_{n, m}(x)=
    (-1)^{n+m} q^{n(n-1)/2} q^{m(m-1)/2} \frac{(q-1)^{n+m}}
    {(q;q)_n (q;q)_m}
    p^{(n+m+1)m}\\
    \times \sum_{k=0}^m
    \qbinom{m}{k}_q  (-1)^k q^{-mk +
    k(k+1)/2}  \frac{(q;q)_{n+m+k}}{(q;q)_{n+k}(1-p)^m}
     D_{q^{-1}}^n[x^{n}
     (qx;q)_{n+k}  ].
\end{multline*}
Using (\ref{artw13}) again, this time with argument $q^{n+k+1}x$
and exponent $n$, gives
\[
    x^n= p^{(n+k+1)n}\sum_{j=0}^n
    \qbinom{n}{j}_q  (-1)^j q^{-nj +
    j(j+1)/2}(q^{n+k+1}x;q)_j,
\]
and since $(qx;q)_{n+k}\, (q^{n+k+1}x;q)_j=(qx;q)_{n+k+j}$ we
obtain
\begin{multline*}
    p_{n, m}(x)=
    (-1)^{n+m} q^{n(n-1)/2} q^{m(m-1)/2} \frac{(q-1)^{n+m}}
    {(q;q)_n (q;q)_m}
    p^{(n+m+1)m}
    \sum_{k=0}^m
    \qbinom{m}{k}_q  (-1)^k q^{-mk +
    k(k+1)/2}\\ \times  \frac{(q;q)_{n+m+k}}{(q;q)_{n+k}(1-p)^m} p^{(n+k+1)n}\sum_{j=0}^n
    \qbinom{n}{j}_q  (-1)^j q^{-nj +
    j(j+1)/2}
     D_{q^{-1}}^n (qx;q)_{n+k+j} .
\end{multline*}
Again using formula (\ref{artw15}) and also using that
\[
    (q;q)_k=(-1)^k p^{-k(k+1)/2}(p;p)_k , \quad
    \qbinom{n}{k}_q= p^{-k(n-k)}
    \qbinom{n}{k}_p,
\]
we can rewrite the expression for $ p_{n, m}(x)$ as
\begin{multline}\label{pnmnewbasis}
    p_{n, m}(x)=
    (-1)^{n+m} \sum_{k=0}^m  \sum_{j=0}^n (-1)^{k+j}
    \qbinom{n+m+k}{m}_p
    \qbinom{n+k+j}{n}_p
    \qbinom{m}{k}_p
    \qbinom{n}{j}_p\\
    \times (qx;q)_{k+j} p^{(n-j)(n-j+1)/2} p^{(m-k)(m-k+1)/2} .
\end{multline}

\section{Hermite-Pad\'{e} approximation of $f_1$ and
$f_2$} \label{HPbenad}

We define two measures $\mu_1$ and $\mu_2$ by taking $d\mu_1(x)=d_qx$
and $d\mu_2(x)=\log_q(x)\,d_qx$, where $d_q$ is defined by (\ref{eq:qint}). 
Then $\mu_1$ and $\mu_2$ are
supported on $\{q^k, k=0,1,2, \ldots\}$, which is a bounded set in
$[0,1]$ with one accumulation point at $0$. The Markov functions
for the measures $\mu_1$ and $\mu_2$ are
\begin{equation}\label{f1}
    f_1(z)= \int_0^1 \frac{d\mu_1(x)}{z-x}= \sum_{k=0}^\infty
    \frac{q^k}{z-q^k}
\end{equation}
\begin{equation}\label{f2}
    f_2(z)= \int_0^1 \frac{d\mu_2(x)}{z-x}= \sum_{k=0}^\infty \frac{k
    q^k}{z-q^k}.
\end{equation}
Observe that for every $N \in \mathbb{N}$
\begin{eqnarray*}
   \zeta_q(1) & = & f_1(p^N) + \sum_{k=1}^{N-1} \frac{1}{p^k-1}, \\
   \zeta_q(2) & = & f_2(p^N) + \sum_{k=1}^{N-1} \frac{k}{p^k-1} + N f_1(p^N),
\end{eqnarray*}
therefore we will look for rational approximants of $f_1(z)$ and $f_2(z)$ with
common denominator and evaluate these at $z=p^N$ for an appropriate choice of $N$.
We can find these rational approximants by Hermite-Pad\'e
approximation of type II.

For Hermite-Pad\'{e} approximation of type II one requires a
polynomial $p_{n,m}$ of degree $\leq n+m$, and
polynomials $q_{n,m}$ and $r_{n,m}$ such that
\begin{equation}\label{Orest1}
    p_{n,m}(z) f_1(z) - q_{n,m}(z)= \mathcal{O}\left(\frac{1}{z^{n+1}}\right), \qquad
    z \to \infty,
\end{equation}
\begin{equation}\label{Orest2}
    p_{n,m}(z) f_2(z) - r_{n,m}(z)= \mathcal{O}\left(\frac{1}{z^{m+1}}\right). \qquad
    z \to \infty,
\end{equation}
It is known \cite[Chapter 4]{NikiSor}
that for $m\leq n$ the polynomial $p_{n,m}$ is, up to a multiplicative factor, uniquely given by
\begin{eqnarray}\label{vwPn}
               \int_0^1 p_{n,m}(x) x^\ell\, d\mu_1(x) & = &0,
                \qquad \ell=0,1, \ldots, n-1, \label{vwPn1} \\
                \int_0^1 p_{n,m}(x) x^\ell\, d\mu_2(x) & = & 0,
                \qquad \ell=0,1, \ldots, m-1,   \label{vwPn2}
 \end{eqnarray}
and that $q_{n,m}$ and $r_{n,m}$ are given by
\begin{equation}\label{QRalgemeen}
    q_{n,m}(z)=
    \int_0^1 \frac{p_{n,m}(z)-p_{n,m}(x)}{z-x}\, d\mu_1(x), \qquad  
    r_{n,m}(z)= 
    \int_0^1 \frac{p_{n,m}(z)-p_{n,m}(x)}{z-x}\, d\mu_2(x).
\end{equation}
The remainder in the approximation (\ref{Orest1}) is
\begin{equation} \label{rest1}
    p_{n,m}(z) f_1(z) -q_{n,m}(z)= \int_0^1 \frac{p_{n,m}(x)}{z-x}\, d\mu_1(x),
\end{equation}
and for (\ref{Orest2})
\begin{equation}\label{rest2}
    p_{n,m}(z) f_2(z) - r_{n,m}(z)=  \int_0^1 \frac{p_{n,m}(x)}{z-x}\, d\mu_2(x).
\end{equation}

Comparing (\ref{vwPn1})--(\ref{vwPn2}) with (\ref{vwpn1})--(\ref{vwpn2}) we see that the common
denominator is given by the multiple little $q$-Jacobi
polynomial $p_{n,m}$. So (\ref{pnmnewbasis}) gives an
explicit expression for $p_{n,m}$. Then it follows that we can compute
$q_{n,m}$ and $r_{n,m}$ explicitly using (\ref{QRalgemeen}). For $q_{n,m}$ we
have to compute
\begin{equation}\label{qnm}
    q_{n,m}(z)=
    \sum_{\ell=0}^\infty
    \frac{p_{n,m}(z)-p_{n,m}(q^\ell)}{z-q^\ell}\ q^\ell.
\end{equation}
By using the explicit expression (\ref{pnmnewbasis}) for $p_{n,m}(z)$, we find
\begin{multline*}
    q_{n, m}(z)=
    (-1)^{n+m} \sum_{k=0}^m  \sum_{j=0}^n (-1)^{k+j}
    \qbinom{n+m+k}{m}_p
    \qbinom{n+k+j}{n}_p
    \qbinom{m}{k}_p
    \qbinom{n}{j}_p\\
    \times  p^{(n-j)(n-j+1)/2} p^{(m-k)(m-k+1)/2}
    \sum_{\ell=0}^\infty
    \frac{(qz;q)_{k+j}-(q^{\ell+1};q)_{k+j}}{z-q^\ell}q^\ell.
\end{multline*}
Now use
\[
    \frac{(qx;q)_{k}-(qy;q)_{k}}{x-y}= -
    \sum_{\ell=1}^{k}q^\ell(qy;q)_{\ell-1}(q^{\ell+1}x;q)_{k-\ell},
\]
which one can prove by induction, then this gives
\begin{multline*}
    q_{n, m}(z)=
    (-1)^{n+m+1} \sum_{k=0}^m  \sum_{j=0}^n (-1)^{k+j}
    \qbinom{n+m+k}{m}_p
    \qbinom{n+k+j}{n}_p
    \qbinom{m}{k}_p
    \qbinom{n}{j}_p\\
    \times  p^{(n-j)(n-j+1)/2} p^{(m-k)(m-k+1)/2}
    \sum_{r=1}^{k+j}q^r (q^{r+1}z;q)_{k-r+j}  \sum_{\ell=0}^\infty
    q^\ell (q^{\ell+1};q)_{r-1}.
\end{multline*}
By using the $q$-binomial series (\ref{artw14}), we can compute
the modified moments
\begin{equation} \label{modifiedmoments}
     \sum_{\ell=0}^\infty
    q^\ell (q^{\ell+1};q)_{r-1}= \frac{1}{1-q^r},
\end{equation}
so that
\begin{multline}\label{qnm uitdrukking}
    q_{n, m}(z)=
    (-1)^{n+m+1} \sum_{k=0}^m  \sum_{j=0}^n (-1)^{k+j}
    \qbinom{n+m+k}{m}_p
    \qbinom{n+k+j}{n}_p
    \qbinom{m}{k}_p
    \qbinom{n}{j}_p\\
    \times  p^{(n-j)(n-j+1)/2} p^{(m-k)(m-k+1)/2}
    \sum_{r=1}^{k+j} \frac{(q^{r+1}z;q)_{k-r+j}}{p^r-1}.
\end{multline}

For an explicit expression of $r_{n,m}$ we use (\ref{QRalgemeen}), which gives
\begin{equation}\label{rnm}
    r_{n,m}(z)=
    \sum_{\ell=0}^\infty 
    \frac{p_{n,m}(z)-p_{n,m}(q^\ell)}{z-q^\ell}\ \ell q^\ell.
\end{equation}
Completely analogous to $q_{n,m}$ we find
\begin{multline*}
    r_{n, m}(z)=
    (-1)^{n+m+1} \sum_{k=0}^m  \sum_{j=0}^n (-1)^{k+j}
    \qbinom{n+m+k}{m}_p
    \qbinom{n+k+j}{n}_p
    \qbinom{m}{k}_p
    \qbinom{n}{j}_p\\
    \times  p^{(n-j)(n-j+1)/2} p^{(m-k)(m-k+1)/2}
    \sum_{r=1}^{k+j}q^r (q^{r+1}z;q)_{k-r+j}  \sum_{\ell=0}^\infty
    \ell q^\ell (q^{\ell+1};q)_{r-1}.
\end{multline*}
Now we have more work to compute
\begin{equation}\label{tussenres}
    \sum_{\ell=0}^\infty \ell q^\ell (q^{\ell+1};q)_{r-1}= (q;q)_{r-1}
    \sum_{\ell=0}^\infty \ell q^\ell \frac{(q^r;q)_\ell}{(q;q)_\ell}.
\end{equation}
Using the $q$-binomial series (\ref{artw14}), we have that
\[
    \quad (x;q)_r \, \sum_{\ell=0}^\infty \frac{(q^r;q)_\ell}{(q;q)_\ell} x^\ell = 1.
\]
By taking the derivative with respect to $x$ we find
\[
    \frac{d}{dx}[(x;q)_r] \sum_{\ell=0}^\infty \frac{(q^r;q)_\ell}{(q;q)_\ell} x^\ell
   + (x;q)_r  \sum_{\ell=0}^\infty \ell \frac{(q^r;q)_\ell}{(q;q)_\ell}  x^{\ell-1} = 0.
\]
It is not difficult to see that  
\[  \frac{d}{dx}[(x;q)_r]= -
   \sum_{i=0}^{r-1}\frac{(x;q)_r}{1-xq^i}\ q^i, \]  
so we find
\[
    \sum_{\ell=0}^\infty \ell \frac{(q^r;q)_\ell}{(q;q)_\ell} x^{\ell}
    = x \sum_{i=0}^{r-1}\frac{q^i}{1-xq^i}
    \sum_{\ell=0}^\infty \frac{(q^r;q)_\ell}{(q;q)_\ell} x^\ell.
\]
From this it follows that (\ref{tussenres}) becomes
\[
    \sum_{\ell=0}^\infty \ell q^\ell (q^{\ell+1};q)_{r-1}= (q;q)_{r-1}
     q \sum_{i=0}^{r-1}\frac{q^i}{1-q^{i+1}}
    \sum_{\ell=0}^\infty \frac{(q^r;q)_\ell}{(q;q)_\ell}\  q^\ell.
\]
Using expression (\ref{modifiedmoments}) for the modified moments
we find
\[
    \sum_{\ell=0}^\infty \ell q^\ell (q^{\ell+1};q)_{r-1}=
    \frac{1}{1-q^r} \sum_{i=1}^{r}\frac{1}{p^{i}-1}.
\]
When we use this in the expression for $r_{n,m}$ we get
\begin{multline}\label{rnm uitdrukking}
    r_{n, m}(z)=
    (-1)^{n+m+1} \sum_{k=0}^m  \sum_{j=0}^n (-1)^{k+j}
    \qbinom{n+m+k}{m}_p
    \qbinom{n+k+j}{n}_p
    \qbinom{m}{k}_p
    \qbinom{n}{j}_p\\
    \times  p^{(n-j)(n-j+1)/2} p^{(m-k)(m-k+1)/2}
    \sum_{r=1}^{k+j} \sum_{i=1}^{r}
    \frac{(q^{r+1}z;q)_{k-r+j}}{(p^r-1)(p^{i}-1)} .
\end{multline}

We will evaluate these functions $p_{n,m}(z)$, $q_{n,m}(z)$,
$r_{n,m}(z)$ at $z=p^{n+m}$. First we will show that $p_{n,m}(p^{n+m})$
is an integer. From the $q$-version of
Pascal's triangle identity
\begin{equation}\label{Pascal}
    \qbinom{n}{k}_p
    =\qbinom{n-1}{k-1}_p + p^k \qbinom{n-1}{k}_p
    =\qbinom{n-1}{k}_p + p^{n-k} \qbinom{n-1}{k-1}_p,
\end{equation}
it follows (by induction) that $\qbinom{n}{k}_p$ is an integer whenever
$p$ is an integer. Furthermore we have
\[
(qp^{n+m};q)_{k+j}=(p^{n+m-1};p^{-1})_{k+j}=\prod_{i=1}^{k+j}(1-p^{n+m-i}).
\]
For each value of $i$ in this product, the factor $1-p^{n+m-i}$ is an integer. This means that
$p_{n,m}(p^{n+m})$, with $p_{n,m}$ given by
(\ref{pnmnewbasis}), is an integer.

Now we will evaluate $q_{n,m}(z)$ and $r_{n,m}(z)$ at $p^{n+m}$.
We can use (\ref{qnm uitdrukking}) at $z=p^{n+m}$ and
    $(q^{r+1}p^{n+m};q)_{k-r+j}=(p^{n+m-k-j};p)_{k-r+j}$, to find
\begin{multline}\label{qnm ge-evalueerd}
    q_{n, m}(p^{n+m})=
    (-1)^{n+m+1} \sum_{k=0}^m  \sum_{j=0}^n (-1)^{k+j}
    \qbinom{n+m+k}{m}_p
    \qbinom{n+k+j}{n}_p
    \qbinom{m}{k}_p
    \qbinom{n}{j}_p\\
    \times  p^{(n-j)(n-j+1)/2} p^{(m-k)(m-k+1)/2}
    \sum_{r=1}^{k+j} \frac{(p^{n+m-k-j};p)_{k-r+j}}{p^r-1}.
\end{multline}
The terms in the sum for $q_{n, m}(p^{n+m})$ are not all integers,
because of the expression $p^r-1$ in the denominators. In order to
obtain an integer we have to multiply $q_{n, m}(p^{n+m})$ by a
multiple of all $p^r-1$ for $r=1, 2, \ldots, n+m$.

For $r_{n, m}(p^{n+m})$ we use (\ref{rnm uitdrukking})
at $z=p^{n+m}$ and again  $(q^{r+1}p^{n+m};q)_{k-r+j}=(p^{n+m-k-j};p)_{k-r+j}$, to find
\begin{multline}\label{rnm ge-evalueerd}
    r_{n,m}(p^{n+m})=
    (-1)^{n+m+1} \sum_{k=0}^m  \sum_{j=0}^n (-1)^{k+j}
    \qbinom{n+m+k}{m}_p
    \qbinom{n+k+j}{n}_p
    \qbinom{m}{k}_p
    \qbinom{n}{j}_p\\
    \times  p^{(n-j)(n-j+1)/2} p^{(m-k)(m-k+1)/2}
    \sum_{r=1}^{k+j} \sum_{i=1}^{r} \frac{(p^{n+m-k-j};p)_{k-r+j}} {(p^r-1)(p^{i}-1)}.
\end{multline}
We see that in order to get an integer now, we have to multiply
$r_{n,m}(p^{n+m})$ by a multiple of $(p^r-1)(p^i-1)$ for $r=1, 2, \ldots,
n+m$ and $i=1, 2, \ldots, n+m$.

Define
\begin{equation}\label{dn}
    d_{n}(x)=\prod_{k=1}^{n} \Phi_k(x), 
\end{equation}
where
\begin{equation} \label{cyclo}
    \Phi_n(x)= \prod_{\substack{k=1\\[0,5ex] \gcd(k,n)=1}}^{n} (x-e^{2\pi i k/n})
\end{equation}
are the cyclotomic polynomials. Each cyclotomic polynomial is
monic and has integer coefficients. It is known \cite{Stillwell} that
\begin{equation}\label{x^n-1}
    x^n-1=\prod_{d|n} \Phi_d(x),
\end{equation}
and that every cyclotomic polynomial is irreducible over
$\mathbb{Q}[x]$. Hence $d_{n}(p)$ is a multiple of all $p^\ell-1$ for
$\ell=1,2, \ldots, n$. The growth of this sequence is given by the
following lemma \cite[Lemma 2]{VanAssche}.

\begin{lemma}\label{dn-asymptotiek}
Suppose $x >1$ and let $d_n(x)$ be
given by (\ref{dn}). Then
\[
    \lim_{n \to \infty}d_n(x)^{1/n^2}=x^{3/\pi^2}.
\]
\end{lemma}

Since $d_{n+m}(p)$ is a multiple of $p^\ell-1$ for all $\ell=1,2,
\ldots, n+m$, we conclude that
     $d_{n+m}(p) q_{n,m}(p^{n+m})$  and
     $d_{n+m}^2(p) r_{n,m}(p^{n+m})$
are integers.

\section{Rational approximations to $\zeta_q(1)$}
 From (\ref{qzetas}) we know that
\begin{equation*}
    \zeta_q(1)=\sum_{n=1}^\infty \frac{q^n}{1-q^n}= \sum_{n=1}^\infty
    \frac{1}{p^n-1}.
\end{equation*}
In this section we will construct rational approximations to
$\zeta_q(1)$, prove the irrationality and give an upper bound for
its measure of irrationality.

\subsection{Rational approximations}
In (\ref{f1}) we defined the function $f_1$. When we evaluate
$f_1$ at $p^{n+m}$, we get
\begin{equation}\label{eq:f1geëvalueerd}
    f_1(p^{n+m})=\sum_{k=0}^\infty
    \frac{q^k}{p^{n+m}-q^k}=\sum_{k=0}^\infty\frac{1}{p^{n+m+k}-1}
    = \zeta_q(1)-\sum_{k=1}^{n+m-1}\frac{1}{p^k-1}
\end{equation}
and hence
\begin{equation}\label{qzeta1}
    \zeta_q(1)=f_1(p^{n+m})+\sum_{k=1}^{n+m-1}\frac{1}{p^k-1}.
\end{equation}
We now take $m=n-1$ and define
\begin{eqnarray}
    \beta_{n}&=& d_{2n-1}(p)p_{n,n-1}(p^{2n-1})\label{beta}, \\
    \alpha_{n}&=&  d_{2n-1}(p)\left[q_{n,n-1}(p^{2n-1})+
    p_{n,n-1}(p^{2n-1})\sum_{k=1}^{2n-2}\frac{1}{p^k-1}\right]\label{alpha}.
\end{eqnarray}
Then it follows from equation (\ref{qzeta1}) and
from the Hermite-Pad\'{e} approximation of $f_1$ (\ref{rest1})
that 
\begin{eqnarray}
    \beta_{n} \zeta_q(1) - \alpha_{n}&=&
    d_{2n-1}(p) [p_{n,n-1}(p^{2n-1})f_1(p^{2n-1})-q_{n,n-1}(p^{2n-1})] \nonumber\\
    &=& d_{2n-1}(p) \int_0^1 \frac{p_{n,n-1}(x)}{p^{2n-1}-x}\ d_qx. \label{betan qzeta1 - alfan}
\end{eqnarray}

\subsection{Irrationality of $\zeta_q(1)$}
Using Lemma \ref{lem:irrational}, we can prove the irrationality
of $\zeta_q(1)$. We will first show that
\[
    \lim_{n \to \infty} |\beta_{n} \zeta_q(1) - \alpha_{n}|=0,
\]
i.e., using (\ref{betan qzeta1 - alfan}), we want to find an
estimate for
\[
    \int_0^1 \frac{p_{n,n-1}(x)}{z-x}\ d_qx.
\]
When we use the Rodrigues formula (\ref{rodriguesa}) for
$p_{n,m}$
 we get
\[
    \int_0^1 \frac{p_{n,m}(x)}{z-x}\ d_qx= C \int_0^1 D_{q^{-1}}^n[x^{n}
    D_{q^{-1}}^m \left((qx;q)_{n+m} x^m \right)]     \frac{1}{z-x}\ d_qx ,
\]
where
 \begin{equation}\label{C}
    C=(-1)^{n+m} q^{n(n-1)/2} q^{m(m-1)/2} \frac{(q-1)^{n+m}} {(q;q)_n
(q;q)_m}.
\end{equation}

Repeated application of summation by parts (\ref{partsomp_q})
gives 
\[
    \int_0^1 \frac{p_{n,m}(x)}{z-x}\ d_qx
    = C (-q)^n \int_0^1 x^{n} D_{q^{-1}}^m \left((qx;q)_{n+m} x^m \right) D_q^n \left[
   \frac{1}{z-x}\right]\, d_qx.
\]
Now it is easy to see by induction that
\[
     D_q^n \left[ \frac{1}{z-x}\right]= \frac{(q;q)_n
     (1-q)^{-n}}{(z-x)(z-xq)\ldots (z-xq^n)}=\frac{(q;q)_n}{
     (1-q)^{n}z^{n+1}} \frac{1}{(\frac{x}{z};q)_{n+1}},
\]
so for $z=p^{n+m}$ we find
\[
    \int_0^1 \frac{p_{n,m}(x)}{p^{n+m}-x}\ d_qx= C (-q)^n
    \frac{(q;q)_n}
     {(1-q)^{n}p^{(n+m)(n+1)}}\int_0^1 \frac{x^{n}
    D_{q^{-1}}^m \left((qx;q)_{n+m} x^m \right)}{(q^{n+m}x;q)_{n+1}}\ d_qx.
\]
Again we apply summation by parts (\ref{partsomp_q}) several times
and use the $q$-analogue of Leibniz' formula
(\ref{Leibniz}). This gives
\begin{multline*}
    \int_0^1 \frac{p_{n,m}(x)}{p^{n+m}-x}\ d_qx= C (-q)^{n+m}
    \frac{(q;q)_n} {(1-q)^{n}p^{(n+m)(n+1)}}\\ 
    \times \sum_{k=0}^m \qbinom{m}{k}_q
         \int_0^1 (qx;q)_{n+m} x^m  \left. \rule{0pt}{15pt} D_q^{m-k}\left(x^n\right) \right|_{q^k x} \,
    D_q^k  \frac{1 }{(q^{n+m}x;q)_{n+1}}\ d_qx.
\end{multline*}
By induction, it is easy to see that for $r \leq s$
\begin{equation}\label{eq:Dqr x^s}
    D_q^{r}\left(x^s\right)= \frac{(q;q)_s \, x^{s-r}}{(1-q)^{r} (q;q)_{s-r}},
\end{equation} 
and that
\[
    D_q^k \frac{1 }{(q^{n+m}x;q)_{n+1}}= \frac{q^{k(n+m)} (q^{n+1};q)_k}
        {(1-q)^k (q^{n+m}x;q)_{n+k+1}}.
\]
Using this we find
\begin{multline}\label{beta qzeta1 -alpha neq 0}
    \int_0^1 \frac{p_{n,m}(x)}{p^{n+m}-x}\ d_qx= C (-q)^{n+m}
    \frac{(q;q)_n} {(1-q)^{n}p^{(n+m)(n+1)}} \\ 
   \times \sum_{k=0}^m \qbinom{m}{k}_q
     \frac{ q^{k(2n+k)}(q;q)_n (q^{n+1};q)_k}{(1-q)^{m} (q;q)_{n-m+k}}
         \int_0^1 \frac{(qx;q)_{n+m}  x^{n+k}}{(q^{n+m}x;q)_{n+k+1}}\ d_qx.
\end{multline}
The integrand is always positive and $(q^{n+m}x;q)_{n+k+1}\geq (q^{n+m};q)_{n+k+1}$, hence we find
\begin{eqnarray*}
    \int_0^1 \frac{(qx;q)_{n+m}  x^{n+k}}{(q^{n+m}x;q)_{n+k+1}}\ d_qx 
     & \leq &  
    \frac{1}{(q^{n+m};q)_{n+k+1}}\int_0^1 (qx;q)_{n+m}  x^{n+k}\, d_qx \\
    & \leq & 
    \frac{1}{(q^{n+m};q)_{n+k+1}}\sum_{\ell=0}^\infty (q^{\ell+1};q)_{n+m}  (q^\ell)^{n+k+1}    \\ 
   & \leq & \frac{(q;q)_{n+m}}{(q^{n+m};q)_{n+k+1}} \sum_{\ell=0}^\infty \frac{(q^{n+m+1};q)_\ell}
          {(q;q)_{\ell}}  (q^{n+k+1})^{\ell}.
\end{eqnarray*}
Using the $q$-binomial series (\ref{artw14}), we then find 
\begin{eqnarray*}
  \int_0^1 \frac{(qx;q)_{n+m}  x^{n+k}}{(q^{n+m}x;q)_{n+k+1}}\ d_qx 
  & \leq &  \frac{(q;q)_{n+m}}{(q^{n+m};q)_{n+k+1}}
         \frac{(q^{2n+k+m+2};q)_\infty}{(q^{n+k+1};q)_\infty} \\
  & = & \frac{(q;q)_{n+k}}{(q^{n+m};q)_{n+k+1}(q^{n+m+1};q)_{n+k+1}} .
\end{eqnarray*}
This estimate, together with the definition of $C$ in (\ref{C}), then
gives 
\begin{multline*}
    \left|\int_0^1 \frac{p_{n,m}(x)}{p^{n+m}-x}\ d_qx\right| \leq 
     \frac{q^{n(n-1)/2}  q^{m(m-1)/2}q^{n+m}} {p^{(n+m)(n+1)}  (q;q)_m}  \\ 
   \times \sum_{k=0}^m \qbinom{m}{k}_q
     \frac{ q^{k(2n+k)}(q;q)_n (q^{n+1};q)_k  (q;q)_{n+k} }
     {(q;q)_{n-m+k} (q^{n+m};q)_{n+k+1}(q^{n+m+1};q)_{n+k+1} } .
\end{multline*}
Some simple estimations then give
\begin{eqnarray*}
    \left|\int_0^1 \frac{p_{n,m}(x)}{p^{n+m}-x} d_qx\right| & \leq &
     \frac{q^{n(n-1)/2}  q^{m(m-1)/2}q^{n+m}} {
     p^{(n+m)(n+1)}  (q;q)_m}  \frac{1} { (q^{n+m};q)_{n+m+2}^2 }
      \sum_{k=0}^m \left(q^{2n}\right)^k \\
     & = & \frac{q^{n(n-1)/2}  q^{m(m-1)/2}q^{n+m}} {
     p^{(n+m)(n+1)}  (q;q)_m  (q^{n+m};q)_{n+m+2}^2 } \frac{1-q^{2n(m+1)}}{1-q^{2n}}.
\end{eqnarray*} 
This gives a useful estimate for the integral on the right hand side of
equation (\ref{betan qzeta1 - alfan}), which (for $m=n-1$) implies that
\begin{equation}\label{beta qzeta1- alpha afschatting}
    \left| \beta_n \zeta_q(1) - \alpha_n \right| <
     d_{2n-1}(p) \,
     \frac{q^{n(n-1)/2}  q^{(n-1)(n-2)/2}} {
     p^{(2n-1)(n+1)} }
    \frac{q^{2n-1}} {(q;q)_{n-1} (q^{2n-1};q)_{2n+1}^2 (1-q^{2n})} .
\end{equation}
Now we can prove Theorem~\ref{thm:qzeta1}

\textbf{Proof of Theorem~\ref{thm:qzeta1}.}
From Section \ref{HPbenad} we know that $\alpha_n$ and $\beta_n$
are integers. From the equations (\ref{betan qzeta1 - alfan}) and
(\ref{beta qzeta1 -alpha neq 0}) it follows that $\beta_n
\zeta_q(1) - \alpha_n \neq 0$ because $d_{2n-1}(p)\neq 0$ and
because the integral on the right hand side of (\ref{betan qzeta1
- alfan}) can be written as a sum with all terms different from
zero and of the same sign. From  (\ref{beta qzeta1- alpha
afschatting}) we can find that
\begin{eqnarray*}
    \lim_{n \to \infty} \left| \beta_n \zeta_q(1) - \alpha_n \right|^{1/n^2}
    &\leq&  \lim_{n \to \infty} d_{2n-1}(p)^{1/n^2} \  q^3 \\
    &  = & p^{12/\pi^2} p^{-3} = p^{-\frac{3(\pi^2-4)}{\pi^2}} <1,
\end{eqnarray*}
where we have used Lemma \ref{dn-asymptotiek}. Hence $\lim_{n \to \infty} \left| \beta_n \zeta_q(1) -
\alpha_n \right| \to 0$. The irrationality now
follows from Lemma \ref{lem:irrational}.  \hspace*{\fill} \qed

\subsection{Measure of irrationality for $\zeta_q(1)$}
Theorem \ref{thm:qzeta1} gives rational approximations
$\alpha_n/\beta_n$ for $\zeta_q(1)$ that satisfy
\begin{equation}\label{qzeta1-alpha/beta}
    \left|\zeta_q(1) - \frac{\alpha_{n}}{\beta_{n}}\right|
    = \mathcal{O} \left( \frac{ p^{\left(\frac{-3(\pi^2-4)}{\pi^2}+\epsilon\right)n^2}}{\beta_{n}}\right)
\end{equation}
for every $\epsilon > 0$, with $\beta_{n}= d_{2n-1}(p)p_{n,n-1}(p^{2n-1})$. 
We already know the asymptotic behavior of
$d_{2n-1}(p)$, so what remains is to find the asymptotic behavior
of $p_{n,n-1}(p^{2n-1})$ as $n \to \infty$.

Since $\{1, x, x^2, \ldots, x^{n-1}, \log x, x\log x, x^2 \log
x, \ldots, x^{m-1} \log x \}$ is a Chebyshev system on $[0,1]$ whenever $m \leq n$, we know that
all the zeros of $p_{n,m}(z)$ are simple and lie in the interval
$(0,1)$. If we call these zeros $x_1, x_2, \ldots, x_{n+m}$ then
we can write $p_{n,m}(x)$ as
\[
    p_{n,m}(x)= \kappa_{n,m} \prod_{j=1}^{n+m}(x-x_j).
\]
For $|x| > 1$ we have $|x|-1 \leq |x-x_j| \leq |x|+1$, hence
\[
    (|x|-1)^{n+m} \leq \prod_{j=1}^{n+m} |x-x_j| \leq
    (|x|+1)^{n+m},
\]
so multiplying by $|\kappa_{n,m}|$ and evaluating at $x=p^{n+m}$
gives
\[
    |\kappa_{n,m}| (p^{n+m}-1)^{n+m} \leq  |p_{n,m}(p^{n+m})|
    \leq |\kappa_{n,m}| (p^{n+m}+1)^{n+m}.
\]
For $m=n-1$ one obtains
\[
    |\kappa_{n,n-1}|^{1/n^2} \left(p^{2n-1}-1\right)^{\frac{2n-1}{n^2}}
    \leq  |p_{n,n-1}(p^{2n-1})|^{1/n^2} \leq
   |\kappa_{n,n-1}|^{1/n^2} \left(p^{2n-1}+1\right)^{\frac{2n-1}{n^2}},
\]
so we have that
\[
    \lim_{n \to \infty} |p_{n,n-1}(p^{2n-1})|^{1/n^2} =
    p^4 \lim_{n \to \infty} |\kappa_{n,n-1}|^{1/n^2}.
\]

Equation (\ref{kappainverse}) implies that
the leading coefficient of $p_{n,m}$ is given by
\begin{equation}\label{Kappa1}
    \kappa_{n,m}=[p_{n, m}^{(0, 0, 0)}(0)]^{-1}
    = (-1)^{n+m}
    \frac    {(q^{n+m+1};q)_n
    (q^{n+m+1};q)_m} { q^{nm} (q;q)_n (q;q)_m q^{n(n-1)/2} q^{m(m-1)/2}}.
\end{equation}
Since $(q^{n+1};q)_k=(q;q)_{n+k}/(q;q)_n$, we can rewrite the
leading coefficient of $p_{n,m}$ as
\begin{equation}\label{Kappa2}
    \kappa_{n,m}
    = (-1)^{n+m}\qbinom{2n+m}{n}_q
    \qbinom{n+2m}{m}_q
     p^{(n+m)^2/2} p^{-(n+m)/2}.
\end{equation}
This gives    $\lim_{n \to \infty} |\kappa_{n,n-1}|^{1/n^2}=
p^2$, so we have
\begin{equation}\label{asymptotiekp}
    \lim_{n \to \infty} |p_{n,n-1}(p^{2n-1})|^{1/n^2} =
    p^{6}.
\end{equation}
Combining this with Lemma \ref{dn-asymptotiek} we have for the
integers $\beta_{n}$ in (\ref{beta}) that
\begin{equation*}
    \lim_{n \to \infty} \left| \beta_n \right|^{1/n^2}=p^{12/\pi^2}p^{6}
    =p^\frac{6(\pi^2+2)}{\pi^2}.
\end{equation*}
Together with (\ref{qzeta1-alpha/beta})  this gives that
\begin{equation}
     \left|\zeta_q(1) - \frac{\alpha_{n}}{\beta_{n}}\right|
    = \mathcal{O} \left( \frac{1}{ \beta_{n}^{1+ \frac{\pi^2-4}{2 (\pi^2+2)}-\epsilon}}\right)
\end{equation}
for every $\epsilon >0$, which implies the following upper bound
for the measure of irrationality
\[
    r(\zeta_q(1)) \leq 1+ \frac{2(\pi^2+2)}{\pi^2-4}= \frac{3
    \pi^2}{\pi^2-4}\approx 5.04443.
\]

\section{Rational approximations to $\zeta_q(2)$}
From (\ref{qzetas}) we know that
\begin{equation*}
    \zeta_q(2)=\sum_{n=1}^\infty \frac{nq^n}{1-q^n}= \sum_{n=1}^\infty
    \frac{n}{p^n-1}.
\end{equation*}
In this section we will construct rational approximations to
$\zeta_q(2)$, we will prove that $\zeta_q(2)$ is irrational, and we
will give an upper bound for its measure of irrationality.

\subsection{Rational approximations}
Evaluating the function $f_2$, which we defined in (\ref{f2}), at
the point $p^{n+m}$, gives 
\[
    f_2(p^{n+m})=\sum_{k=0}^\infty
    \frac{k
    q^k}{p^{n+m}-q^k}=\sum_{k=0}^\infty\frac{k}{p^{n+m+k}-1}.
\]
So we have
\begin{eqnarray*}
    f_2(p^{n+m})
    &=&\sum_{k=0}^\infty\frac{n+m+k}{p^{n+m+k}-1}-(n+m)\sum_{k=0}^\infty\frac{1}{p^{n+m+k}-1}\\
    &=&\zeta_q(2)-\sum_{k=1}^{n+m-1}\frac{k}{p^{k}-1}-(n+m)f_1(p^{n+m}),
\end{eqnarray*}
where we used (\ref{eq:f1geëvalueerd}) for $f_1(p^{n+m})$.
Hence we can write $\zeta_q(2)$ as follows:
\begin{equation}\label{qzeta2}
    \zeta_q(2)=f_2(p^{n+m})+\sum_{k=1}^{n+m-1}\frac{k}{p^{k}-1}+(n+m)f_1(p^{n+m}).
\end{equation}
We now take $m=n-1$ and define
\begin{eqnarray}
  \hspace*{-15pt}  b_{n}&=& d_{2n-1}^2(p)p_{n,n-1}(p^{2n-1}),\label{bnm}\\
   \hspace*{-15pt} a_{n}&=&
    d_{2n-1}^2(p)\left[r_{n,n-1}(p^{2n-1})+
    p_{n,n-1}(p^{2n-1})\sum_{k=1}^{2n-2}\frac{k}{p^k-1} 
      + (2n-1)q_{n,n-1}(p^{2n-1})\right].  \label{anm}
\end{eqnarray}
Then it follows from (\ref{qzeta2}) and
from the Hermite-Pad\'{e} approximation of $f_1$ and $f_2$
(\ref{rest1})--(\ref{rest2}) that
\begin{equation}\label{bz-a_voorlopig}
    b_n \zeta_q(2) - a_n=
    d_{2n-1}^2(p)\left[\int\frac{p_{n,n-1}(x)}{p^{2n-1}-x}\, d\mu_2(x) +
    (2n-1)\int\frac{p_{n,n-1}(x)}{p^{2n-1}-x}\,d\mu_1(x)\right].
\end{equation}
From Section \ref{HPbenad} we know that the numbers $a_n$ and
$b_n$ are integers for all $n \in \mathbb{N}$. We will show that
\[
    \lim_{n \to \infty} |b_n \zeta_q(2) - a_n|=0,
\]
and that $b_n \zeta_q(2) - a_n \neq 0$ for all $n \in \mathbb{N}$,
so that Lemma \ref{lem:irrational} implies the irrationality of
$\zeta_q(2)$. But before we can do that, we have to prove some
results about the asymptotic behavior of $p_{n,m}$ as $n,m \to \infty$.

\subsection{Asymptotic behavior of $p_{n,m}$}
It is well known that the common denominator of Hermite-Pad\'e
approximants satisfies a multiple orthogonality relation. Driver and
Stahl \cite{Driver-Stahl} showed that for a Nikishin system these
common denominators also satisfy an ordinary orthogonality
relation. Although in our case $(f_1, f_2)$ do not form a Nikishin
system, we can prove in a similar way
that $p_{n,m}$ in our case also satisfies an ordinary
orthogonality relation. We need the following theorem.
\begin{theorem}\label{ortheigrest1}
Let $p_{n,m}$ be the multiple little $q$-Jacobi polynomial given
by (\ref{vwpn1})--(\ref{vwpn2}), let $q_{n,m}$ be defined by (\ref{qnm}) and let
$m \leq n$. Then we have
\begin{equation}\label{orthorest}
    \int_{- \infty}^0 (y-1)^k[p_{n,m}(y)f_1(y)-q_{n,m}(y)]\,dy=0
\end{equation}
for $k=0, 1, \ldots, m-1$.
\end{theorem}
\begin{proof}
From the expression of the remainder of the Hermite-Pad\'{e}
approximation (\ref{rest1}) it follows that
\[
    \int_{- \infty}^0 (y-1)^k[p_{n,m}(y)f_1(y)-q_{n,m}(y)]\,dy=
    \int_{- \infty}^0 \frac{(y-1)^{k+1}}{y-1}
    \left[\int_0^1 \frac{p_{n,m}(x)}{y-x} \, d\mu_1(x) \right]dy.
\]
If we add and subtract $(x-1)^{k+1}$ on the right hand side, then
we have
\begin{eqnarray*}
\lefteqn{\int_{- \infty}^0 (y-1)^k[p_{n,m}(y)f_1(y)-q_{n,m}(y)]\, dy} \rule{1in}{0in} & & \\ 
 & = &  \int_{- \infty}^0 \frac{1}{y-1}
    \left[\int_0^1 \frac{\left( (y-1)^{k+1}-(x-1)^{k+1}  \right)}{y-x}
    p_{n,m}(x)\,  d\mu_1(x) \right]dy \\
  & &  +  \int_{- \infty}^0 \frac{1}{y-1}
    \left[\int_0^1 \frac{(x-1)^{k+1}}{y-x} p_{n,m}(x)\, d\mu_1(x) \right]dy.
\end{eqnarray*}
Now $\left( (y-1)^{k+1}-(x-1)^{k+1}  \right)/( y-x)$ is a polynomial
of degree $k$ in $x$, so because of the orthogonality (\ref{vwpn1}) the first integral
on the right hand side vanishes for $k<n$. Since $m \leq n$,
this integral therefore vanishes for $k=0, 1, \ldots, m-1$. Changing the
order of integration gives 
\[
    \int_{- \infty}^0 (y-1)^k[p_{n,m}(y)f_1(y)-q_{n,m}(y)]\,dy=
     \int_0^1 (x-1)^{k} p_{n,m}(x) \int_{- \infty}^0   \frac{(x-1)}{(y-1)(y-x)}\, dy\, d\mu_1(x).
\]
By a partial fraction decomposition we obtain
\[
    \frac{1}{1-y} - \frac{1}{x-y} = \frac{(x-1)}{(y-1)(y-x)},
\]
so we find that
\[
    \int_{- \infty}^0   \frac{(x-1)}{(y-1)(y-x)}\, dy =
    \lim_{t \to
    -\infty} \left[ \int_t^0 \frac{dy}{1-y} - \int_t^0 \frac{dy}{x-y}
    \right] =  \log x.
\]
Using this in the previous expression, we find
\[
    \int_{- \infty}^0 (y-1)^k[p_{n,m}(y)f_1(y)-q_{n,m}(y)]\,dy=
     \int_0^1  (x-1)^{k} p_{n,m}(x) \log x \ d\mu_1(x).
\]
The orthogonality relations (\ref{vwpn2}) then imply that
the right hand side is equal to zero for $k<m$.
\end{proof}

From this theorem it follows that $p_{n,m}(y)f_1(y)-q_{n,m}(y)$
has at least $m$ sign changes on the interval $(-\infty,0)$. 
For suppose $p_{n,m}(y)f_1(y)-q_{n,m}(y)$ has only $r<m$
sign changes on $(-\infty, 0)$, say at $s_1, s_2, \ldots, s_r$,
then we have with $\pi_r(z)=(z-s_1)(z-s_2) \ldots (z-s_r)$ that
\[
    \int_{-\infty}^{0} \pi_r(y) \,
    [p_{n,m}(y)f_1(y)-q_{n,m}(y)]\, dy \neq 0
\]
because the integrand has no sign changes on $(-\infty,0)$. This
gives a contradiction with the fact that this integral can be
written as a linear combination of integrals of the form
\[
    \int_{- \infty}^0 (y-1)^k[p_{n,m}(y)f_1(y)-q_{n,m}(y)]\, dy
\]
with $k \leq r <m$, because from Theorem \ref{ortheigrest1} we
know that these integrals are all zero. So
$p_{n,m}(y)f_1(y)-q_{n,m}(y)$ has at least $m$ sign changes on $(-\infty, 0)$.
The condition $m \leq n$ and (\ref{Orest1}) guarantees that these integrals are finite.

Let $s_1, s_2, \ldots, s_m$ be $m$ points where
$p_{n,m}(z)f_1(z)-q_{n,m}(z)$ changes sign on $(-\infty, 0)$. Then we define the polynomial
$w$ by
\begin{equation} \label{w(z)}
    w(z)=(z-s_1)(z-s_2)\ldots(z-s_m).
\end{equation}
Now we can prove that $p_{n,m}$ also satisfies an ordinary orhogonality relation.

\begin{theorem}\label{newortho}
Let $w$ be defined by (\ref{w(z)}) and let $p_{n,m}$ the
multiple little $q$-Jacobi polynomial defined by (\ref{vwpn1})--(\ref{vwpn2}).
Then
\[
    \int_0^1 p_{n,m}(y) y^k\, \frac{d\mu_1(y)}{w(y)}=0,
\]
for $k=0,1, \ldots, n+m-1$.
\end{theorem}
\begin{proof}
Let $\gamma$ be a closed positively oriented path of integration with
winding number 1 for all its interior points, such that the
interval $[0,1]$ is in the interior of $\gamma$ and the zeros of
$w$ are outside $\gamma$. Using (\ref{rest1}) we have that
\[
     \frac{1}{2 \pi i} \int_\gamma \frac{y^k}{w(y)} \int_0^1 \frac{p_{n,m}(x)}{x-y}\, d\mu_1(x)\, dy
     = \frac{1}{2 \pi i} \int_\gamma y^k [p_{n,m}(y) f_1(y) -q_{n,m}(y)]\, \frac{dy}{w(y)}.
\]
Changing the order of integration on the left hand side gives
\[
    \int_0^1 p_{n,m}(x) \frac{1}{2 \pi i} \int_\gamma \frac{y^k}{x-y}\frac{dy}{w(y)}\, d\mu_1(x)
 = \frac{1}{2 \pi i} \int_\gamma y^k [p_{n,m}(y) f_1(y) -q_{n,m}(y)]\, \frac{dy}{w(y)}.
\]
The only singularity inside $\gamma$ is $x$, hence Cauchy's
formula on the left hand side gives
\[
    \int_0^1  p_{n,m}(x) \frac{x^k}{w(x)}\, d\mu_1(x) =
     \frac{1}{2 \pi i} \int_\gamma y^k [p_{n,m}(y) f_1(y) -q_{n,m}(y)]\, \frac{dy}{w(y)}.
\]
Since all zeros of $w$
are also zeros of $p_{n,m}(y) f_1(y) -q_{n,m}(y)$, the function
$y^k [p_{n,m}(y) f_1(y) -q_{n,m}(y)]/w(y)$ is analytic in
$\mathbb{C}\setminus [0,1]$. Furthermore $\gamma$ encloses the interval
$[0,1]$, hence for $R$ sufficient large 
\[
     \left| \int_0^1 p_{n,m}(x) x^k\, \frac{d\mu_1(x)}{w(x)}\right|
     = \left| \frac{1}{2 \pi i} \int_{\Gamma_R} y^k [p_{n,m}(y) f_1(y)
     -q_{n,m}(y)]\, \frac{dy}{w(y)}\right|,
\]
with $\Gamma_R$ the circle with center 0 and radius $R$. We can
estimate the expression by 
\[
     \left| \int_0^1 p_{n,m}(x) x^k\, \frac{d\mu_1(x)}{w(x)}\right|
     \leq  R \, \max_{y \in \Gamma_R} \left|y^k \frac{p_{n,m}(y) f_1(y)
     -q_{n,m}(y)}{w(y)}\right|.
\]
Using (\ref{Orest1}), we know that the function $y^k$
$[p_{n,m}(y) f_1(y) -q_{n,m}(y)]/w(y)$ has a zero of order
$n+m+1-k$ at infinity. So when $R$ tends to infinity we find
\[
     \left| \int_0^1 p_{n,m}(x) x^k \,\frac{d\mu_1(x)}{w(x)}\right|
     =  \mathcal{O}\left(\frac{1}{R^{n+m-k}}\right).
\]
Therefore
\[
      \int_0^1 p_{n,m}(x) x^k\, \frac{d\mu_1(x)}{w(x)} = 0
\]
when $k \leq n+m-1$.
\end{proof}

Now we have proved that $p_{n,m}$ 
satisfies an ordinary orthogonality relation. Using this,
we can easily show that $p_{n,m}(z) f_1(z) -q_{n,m}(z)$
has exactly $m$ sign changes in the interval $(-\infty, 0)$, namely at the
zeros of $w$. Using (\ref{rest1}) we have
\[
    p_{n,m}(z) f_1(z) -q_{n,m}(z)= \int_0^1
    \frac{p_{n,m}(y)}{z-y}\frac{w(y)}{w(y)}\,d\mu_1(y).
\]
If we add and subtract $w(z)$ on the right hand side, then we get
\[
    p_{n,m}(z) f_1(z) -q_{n,m}(z)=
    \int_0^1\frac{p_{n,m}(y)}{z-y}\frac{w(y)-w(z)}{w(y)}\,d\mu_1(y)
    + w(z) \int_0^1\frac{p_{n,m}(y)}{z-y}\,\frac{d\mu_1(y)}{w(y)}.
\]
The first integral on the right hand side vanishes because of the
orthogonality (Theorem \ref{newortho}), so that
\[
    p_{n,m}(z) f_1(z) -q_{n,m}(z)=
    w(z) \frac{p_{n,m}(z)}{p_{n,m}(z)} \int_0^1\frac{p_{n,m}(y)}{z-y}\,\frac{d\mu_1(y)}{w(y)}.
\]
Now we add and subtract $p_{n,m}(y)$ on the right hand side, so we
get
\begin{multline*}
    p_{n,m}(z) f_1(z) -q_{n,m}(z)=\\
    \frac{w(z)}{p_{n,m}(z)}
    \int_0^1\frac{p_{n,m}(y)}{z-y}(p_{n,m}(z)-p_{n,m}(y))\,
    \frac{d\mu_1(y)}{w(y)}
    +\frac{w(z)}{p_{n,m}(z)} \int_0^1 \frac{p_{n,m}^2(y)}{z-y}\,\frac{d\mu_1(y)}{w(y)}.
\end{multline*}
Again, the first expression on the right hand side vanishes
because of the orthogonality (Theorem
\ref{newortho}), so we find
\begin{equation}\label{uitdrfout1}
    p_{n,m}(z) f_1(z) -q_{n,m}(z)=
    \frac{w(z)}{p_{n,m}(z)} \int_0^1\frac{p_{n,m}^2(y)}{z-y}\,\frac{d\mu_1(y)}{w(y)}.
\end{equation}
The integrand is not identically zero and has a constant sign for $z \in \mathbb{R} \setminus [0,1]$,
hence the integral can not be zero, so it follows from this expression
that the sign changes of $p_{n,m}(z) f_1(z) -q_{n,m}(z)$ in
$\mathbb{R}\setminus [0,1]$ are at the $m$ simple zeros of $w$ in
$(-\infty,0)$.

\subsection{Irrationality of $\zeta_q(2)$}
We will now show that $\lim_{n \to \infty} (b_n \zeta_q(2) -
a_n)=0$ and $b_n \zeta_q(2) - a_n \neq 0$ for all $n \in
\mathbb{N}$, so that Lemma \ref{lem:irrational} implies the
irrationality of $\zeta_q(2)$. First we can reduce the right hand
side of (\ref{bz-a_voorlopig}) as follows. When we look at the
Cauchy transform of  $p_{n,m}(z) f_1(z) -q_{n,m}(z)$ and use
(\ref{rest1}), then we find
\[
\int_{-\infty}^0 \frac{p_{n,m}(z) f_1(z) -q_{n,m}(z)}{y-z}\,dz=
\int_{-\infty}^0 \frac{dz}{y-z} \int_0^1 \frac{p_{n,m}(x)}{z-x}\, d\mu_1(x).
\]
Because
\[
    \frac{1}{(y-z)(z-x)}= \frac{1}{(y-x)}
    \left[\frac{1}{y-z}-\frac{1}{x-z}\right],
\]
 we find by changing the order of integration
\begin{eqnarray}
  \hspace*{-20pt}  \int_{-\infty}^0 \frac{p_{n,m}(z) f_1(z) -q_{n,m}(z)}{y-z}\,dz &=&
    \int_0^1 \frac{p_{n,m}(x)}{y-x}\,d\mu_1(x)    \int_{-\infty}^0 
    \left[\frac{1}{y-z}-\frac{1}{x-z}\right]\, dz \nonumber \\
    &=&
    \int_0^1 \frac{p_{n,m}(x)}{y-x}\log x \ d\mu_1(x) -  \log y \int_0^1 \frac{p_{n,m}(x)}{y-x}\,d\mu_1(x).  
    \label{uitdrmetlog}
\end{eqnarray}
By multiplying both sides by $d_{n+m}^2(p)/\log q$ and using
$\log x \, d\mu_1(x)/ \log q= \log_q x \ d\mu_1(x)= d\mu_2(x)$,
we obtain
\begin{multline*}
    \frac{d_{n+m}^2(p)}{\log q}\int_{-\infty}^0 \frac{p_{n,m}(z) f_1(z) -
    q_{n,m}(z)}{y-z}\, dz\\
    = d_{n+m}^2(p) \left[\int_0^1 \frac{p_{n,m}(x)}{y-x}\, d\mu_2(x)
    - \frac{\log y}{\log  q} \int_0^1 \frac{p_{n,m}(x)}{y-x}\,d\mu_1(x)\right].
\end{multline*}
Evaluating at $z=p^{n+m}$ gives the right hand side of
(\ref{bz-a_voorlopig}) when $m=n-1$, so it turns out that
\begin{equation}\label{bz-a_voorlopig2}
    b_n \zeta_q(2) - a_n=
     \frac{d_{2n-1}^2(p)}{\log q}\int_{-\infty}^0
    \frac{p_{n,n-1}(z) f_1(z) -q_{n,n-1}(z)}{p^{2n-1}-z}\,dz    .
\end{equation}
Now we will estimate the expression on the right hand side.
Multiplying and dividing by $w(p^{n+m})$ and adding and
subtracting $w(z)$ gives that the right hand side is equal to
\begin{multline*}
    \frac{d_{2n-1}^2(p)}{\log q}
    \frac{1}{w(p^{2n-1})} \left[ \int_{-\infty}^0 (w(p^{2n-1})-w(z))
    \frac{p_{n,n-1}(z) f_1(z) -q_{n,n-1}(z)}
    {p^{2n-1}-z}\,dz \right.\\ \left.
    + \int_{-\infty}^0 w(z) \frac{p_{n,n-1}(z) f_1(z) -q_{n,n-1}(z)}{p^{2n-1}-z}\,dz \right].
\end{multline*}
The first term is 0 because of the orthogonality relations given
by Theorem \ref{ortheigrest1}, so we find
\begin{equation} \label{bz-a_nietnul}
    b_n \zeta_q(2) - a_n  =
     \frac{d_{2n-1}^2(p)}{\log q} \frac{1}{w(p^{2n-1})}
       \int_{-\infty}^0 w(z)\frac{p_{n,n-1}(z) f_1(z) -q_{n,n-1}(z)}{p^{2n-1}-z}\,dz.
\end{equation}
Since $w$ is a monic polynomial of degree $m=n-1$ and all the
zeros of $w$ are in the interval $(-\infty, 0)$, we have that
$w(p^{2n-1})> (p^{2n-1})^{n-1}$. We also have that $p^{2n-1}-z >p^{2n-1}$ for $z \in
(-\infty, 0)$, so we find
\begin{equation}\label{bz-a_voorlopig3}
    \left| b_n \zeta_q(2) - a_n \right| \leq -
     \frac{d_{2n-1}^2(p)}{\log q} \left(q^{2n-1}\right)^{n}
      \left| \int_{-\infty}^0 w(z) \left[p_{n,n-1}(z) f_1(z) -q_{n,n-1}(z)\right]\,dz \right| .
\end{equation}
The integral on the right hand side can be evaluated exactly, which is done in the following lemma.
\begin{lemma} \label{extralemma}
Suppose that $m \leq n-1$, then
\begin{equation}\label{integraal w(z) fout1}
    \int_{-\infty}^0 w(z)[p_{n,m}(z) f_1(z) -q_{n,m}(z)] dz=
    q^{n(m+1)} (q;q)_m
     \frac{(q;q)_{n-m-1}}{(q^{n+m+1};q)_{m+1}}
     \frac{(q;q)_m}{(q;q)_{n}}  \log q.
\end{equation}
\end{lemma}
\begin{proof}
Multiply both sides of equation (\ref{uitdrmetlog}) by $w(y)$,
then add and subtract $w(z)$ on both sides of
this expression and use the orthogonality relations
(\ref{orthorest}) for the left hand side and (\ref{vwpn1})--(\ref{vwpn2}) for the
right hand side, to find
\begin{multline*}
    \int_{-\infty}^0 \frac{w(z)}{y-z}[p_{n,m}(z) f_1(z) -q_{n,m}(z)]\,dz=\\
    \int_0^1 \frac{w(z)}{y-z} p_{n,m}(z)\log z \, d\mu_1(z) -
    \log y \int_0^1 \frac{w(z)}{y-z}p_{n,m}(z)d\mu_1(z).
\end{multline*}
Multiply both sides by $y$ and add and subtract $z$ in
the second term on the right hand side to find
\begin{multline*}
    \int_{-\infty}^0 \frac{y}{y-z}w(z)[p_{n,m}(z) f_1(z) -q_{n,m}(z)]\, dz=
    \int_0^1 \frac{y}{y-z}w(z) p_{n,m}(z)\log z \ d\mu_1(z)\\ 
   - \log y \int_0^1 w(z) p_{n,m}(z)\,d\mu_1(z)-
    \log y \int_0^1 \frac{z}{y-z}w(z) p_{n,m}(z)\,d\mu_1(z).
\end{multline*}
We know that $w$ is a polynomial of degree $m$. For
$m\leq n-1$ the orthogonality relations
(\ref{vwpn1})--(\ref{vwpn2}) imply that the second integral on the right hand side is 0.
Taking the limit for $y$ going to $\infty$, we get
\begin{equation*}
    \int_{-\infty}^0 w(z)[p_{n,m}(z) f_1(z) -q_{n,m}(z)]\, dz=
     \int_0^1 w(z) p_{n,m}(z)\log z \ d\mu_1(z).
\end{equation*}
Since $w$ is a monic polynomial of degree $m$, it follows
from the orthogonality relations (\ref{vwpn2}) that  
\begin{equation} \label{nog verder uitrekenen}
    \int_{-\infty}^0 w(z)[p_{n,m}(z) f_1(z) -q_{n,m}(z)]\, dz=
      \int_0^1 z^m p_{n,m}(z) \log z  d\mu_1(z).
\end{equation}
Using the Rodrigues formula (\ref{rodriguesa}) for
$p_{n,m}$ on the right hand side we get
\[
    \int_{-\infty}^0 w(z)[p_{n,m}(z) f_1(z) -q_{n,m}(z)]\, dz 
  = C  \int_0^1 z^m \log z D_{q^{-1}}^n[z^{n}
    D_{q^{-1}}^m \left((qz;q)_{n+m} z^m \right)]\, d\mu_1(z),
\]
with $C$ given by expression (\ref{C}). Repeated application of
summation by parts (\ref{partsomp_q}), gives
\begin{multline*}
    \int_{-\infty}^0 w(z)[p_{n,m}(z) f_1(z) -q_{n,m}(z)]\, dz
  \\ = (-1)^n q^n  C   \int_0^1 D_{q}^n
    \left[z^m \log z\right] z^{n}  D_{q^{-1}}^m \left[(qz;q)_{n+m} {z}^m \right]\, d\mu_1(z).
\end{multline*}
For $D_{q}^n \left[z^m \log z\right]$ we use the $q$-analogue of Leibniz'
    rule (\ref{Leibniz}) to find
\begin{equation}
    D_{q}^n \left[z^m \log z\right] = \sum_{k=0}^n
    \qbinom{n}{k}_q \left. \left( D_{q}^{n-k}
    x^m  \right)\right|_{x=q^{k}z}  \left. \left(D_{q}^k \log
    x\right)\right|_{x=z}.
\end{equation}
Equation (\ref{eq:Dqr x^s}) gives us an expression for
$D_{q}^{n-k} z^m$ when
 $n-k \leq m$. For $n-k>m$ we know that $D_{q}^{n-k}
z^m$ is $0$. It is also easy to see that for $k>0$
\[
    D_{q}^k \log x=(-1)^k q^{-k(k-1)/2} \frac{\log q}{(1-q)^k}
    \frac{(q;q)_{k-1}}{x^k},
\]
so it turns out that
\begin{multline} \label{tussenresultaat som-integraal}
    \int_{-\infty}^0 w(z)[p_{n,m}(z) f_1(z) -q_{n,m}(z)]\, dz
  = (-1)^n q^n  C  \sum_{k=n-m}^n
    \qbinom{n}{k}_q
    \frac{(q;q)_m }{(1-q)^{n}}\frac{\left(q^k
    \right)^{m-n+k}}{(q;q)_{m-n+k}}\\ 
   \times  (-1)^k q^{-k(k-1)/2} (q;q)_{k-1}
      \log q  \int_0^1 z^{m}  D_{q^{-1}}^m \left[(qz;q)_{n+m} {z}^m \right]\, d\mu_1(z).
\end{multline}
We now have to evaluate a $q$-integral and
a finite sum.
\begin{enumerate}
\item  We start with the integral in (\ref{tussenresultaat som-integraal}). Apply summation by parts $m$ times to find
\begin{equation*}
    \int_0^1 z^{m}  D_{q^{-1}}^m \left[(qz;q)_{n+m} {z}^m \right]\, d\mu_1(z)
    = (-1)^m q^m \int_0^1 D_{q}^m [z^{m}]  (qz;q)_{n+m} {z}^m\, d\mu_1(z).
\end{equation*}
When we use equation (\ref{eq:Dqr x^s}) for  $D_{q}^m
[z^{m}]$ and rewrite the $q$-integral as a infinite sum using
(\ref{eq:qint}), it turns out that
\[ 
     \int_0^1 z^{m}
    D_{q^{-1}}^m \left[(qz;q)_{n+m} {z}^m \right]\, d\mu_1(z)=
    (-1)^m q^m \frac{(q;q)_m}{(1-q)^m} \sum_{j=0}^\infty
     \left(q^j\right)^{m+1}(q^{j+1};q)_{n+m}.
\]
We can compute the sum on the right hand side
by using the $q$-binomial series (\ref{artw14}). This gives
\begin{eqnarray*}
     \sum_{j=0}^\infty  \left(q^{m+1}\right)^j (q^{j+1};q)_{n+m}
     &= &(q;q)_{n+m}  \sum_{j=0}^\infty  \left(q^{m+1}\right)^j
     \frac{(q^{n+m+1};q)_{j}}{(q;q)_j}\\
     &= &(q;q)_{n+m}
     \frac{(q^{n+2m+2};q)_{\infty}}{(q^{m+1};q)_\infty}\\
     &=& \frac{(q;q)_m}{(q^{n+m+1};q)_{m+1}}.
\end{eqnarray*}
So we find that
\begin{equation}\label{berekening-integraal}
     \int_0^1 z^{m}  D_{q^{-1}}^m \left[(qz;q)_{n+m} {z}^m \right]\, d\mu_1(z)
    = (-1)^m q^m \frac{(q;q)_m}{(1-q)^m} \frac{(q;q)_m}{(q^{n+m+1};q)_{m+1}}.
\end{equation}
\item Now we will work out the finite sum 
\[
    \sum_{k=n-m}^n
    \qbinom{n}{k}_q
    \frac{(q;q)_m }{(1-q)^{n}}\frac{\left(q^k \right)^{m-n+k}}{(q;q)_{m-n+k}}
    (-1)^k q^{-k(k-1)/2} (q;q)_{k-1}
\]
in (\ref{tussenresultaat som-integraal}).
Changing the index of summation gives
\[
   \sum_{k=0}^m
    \qbinom{n}{k+n-m}_q
    \frac{(q;q)_m }{(1-q)^{n}}\frac{\left(q^{k+n-m} \right)^{k}}{(q;q)_{k}}
    (-1)^{k+n-m} q^{-(k+n-m)(k+n-m-1)/2} (q;q)_{k+n-m-1}.
\]
Writing out the $q$-binomial coefficient and rewriting the powers of $q$  gives 
that this is equal to
\[
   \sum_{k=0}^m
   (-1)^{k+n-m}  \frac{(q;q)_n (q;q)_m  (q;q)_{k+n-m-1}}{(q;q)_{k+n-m}(q;q)_{m-k} (q;q)_{k}}
    \frac{ q^{k(k+1)/2} q^{mn} q^{-n(n-1)/2} q^{-m(m+1)/2}}{(1-q)^{n}}.
\]
Use $(q;q)_{k+n-m-1}/(q;q)_{k+n-m} =1/(1-q^{k+n-m})$ to rewrite this sum as
\[
    (-1)^{n+m} q^{mn} q^{-n(n-1)/2} q^{-m(m+1)/2} \frac{(q;q)_n}{(1-q)^{n}}
   \sum_{k=0}^m  (-1)^k \qbinom{m}{k}_q
   \frac{ q^{k(k+1)/2} }{(1-q^{k+n-m})}.
\]
Now use 
\[  \frac{1}{1-q^{k+n-m}} = \sum_{\ell=0}^\infty \left( q^{k+n-m} \right)^\ell \]
to find
\[
    (-1)^{n+m} q^{mn} q^{-n(n-1)/2} q^{-m(m+1)/2} \frac{(q;q)_n}{(1-q)^{n}}
    \sum_{\ell=0}^\infty \left( q^{n-m} \right)^\ell
    \sum_{k=0}^m  (-1)^k \qbinom{m}{k}_q
    \left( q^{k} \right)^{\ell+1}   q^{k(k-1)/2}.
\]
If we now use the $q$-analog of Newton's binomial formula (\ref{artw12}) with $x=q^{\ell+1}$
then the expression reduces to
\begin{equation*}
    (-1)^{n+m} q^{mn} q^{-n(n-1)/2} q^{-m(m+1)/2} \frac{(q;q)_n}{(1-q)^{n}}
    \sum_{\ell=0}^\infty \left( q^{n-m} \right)^\ell (q^{\ell+1};q)_m .
\end{equation*}
 Using the $q$-binomial series (\ref{artw14}), we have
\begin{eqnarray*}
     \sum_{\ell=0}^\infty  \left(q^{n-m}\right)^\ell (q^{\ell+1};q)_{m}
     & = & (q;q)_{m}  \sum_{\ell=0}^\infty  \left(q^{n-m}\right)^\ell
     \frac{(q^{m+1};q)_{\ell}}{(q;q)_\ell}\\
    & =&(q;q)_{m}
    \frac{(q^{n+1};q)_{\infty}}{(q^{n-m};q)_\infty}\\
    &=& (q;q)_{n-m-1} \frac{(q;q)_{m}}{(q;q)_n}.
\end{eqnarray*}
With this we then find
\begin{multline}
    \sum_{k=n-m}^n
    \qbinom{n}{k}_q
    \frac{(q;q)_m }{(1-q)^{n}}\frac{\left(q^k \right)^{m-n+k}}{(q;q)_{m-n+k}}
    (-1)^k q^{-k(k-1)/2} (q;q)_{k-1} \\
    =  (-1)^{n+m} q^{mn} q^{-n(n-1)/2} q^{-m(m+1)/2} \frac{(q;q)_{m}}{(1-q)^{n}}
     (q;q)_{n-m-1} .
\end{multline}
\end{enumerate}
Using this result together with (\ref{berekening-integraal}) and
the expression  (\ref{C}) for $C$,  we find that
(\ref{tussenresultaat som-integraal}) gives the required result.
\end{proof}

The integral in (\ref{integraal w(z) fout1}) can also be evaluated for $m=n$ but with
much more effort, which is the main reason why we have chosen $m=n-1$ for our rational
approximants. 
Lemma \ref{extralemma} and (\ref{bz-a_voorlopig3}) now imply that for $m=n-1$
\begin{equation}\label{bz-a}
    \left| b_n \zeta_q(2) - a_n \right| \leq
     d_{2n-1}^2(p) \left(q^{2n-1}\right)^{n}
    q^{n^2} 
     \frac{(q;q)_{n-1}^2}{(q^{2n};q)_{n}(q;q)_{n}}.
\end{equation}
Now we can prove Theorem \ref{thm:qzeta2}

{\bf Proof of Theorem \ref{thm:qzeta2}.}
From Section \ref{HPbenad} we know that $a_n$ and $b_n$ are
integers. From  (\ref{bz-a_nietnul}) it follows that $b_n
\zeta_q(2) - a_n \neq 0$ because the integrand has a constant
sign. Furthermore Lemma \ref{dn-asymptotiek} gives that
\[
    \lim_{n \to \infty} \left(d_{2n-1}^2(p)\right)^{1/n^2} = p^{24/\pi^2}.
\]
 From(\ref{bz-a}) we see that
\[
     \lim_{n \to \infty} \left| b_n \zeta_q(2) - a_n \right|^{1/n^2}
      \leq  
     p^\frac{24}{\pi^2}  \lim_{n \to \infty} q^{\frac{2n-1}{n}} q
        = p^\frac{24}{\pi^2}q^3 = p^{-\frac{3(\pi^2-8)}{\pi^2}}<1
\]
so that $\lim_{n \to \infty} \left| b_n \zeta_q(2) - a_n \right|
\to 0$ for $n \to \infty$. The irrationality now follows from
Lemma \ref{lem:irrational}. \hspace*{\fill}$\Box$

\subsection{Measure of irrationality for $\zeta_q(2)$} 
Theorem \ref{thm:qzeta2} gives rational approximations $a_n/b_n$ for
$\zeta_q(2)$ that satisfy
\begin{equation}\label{irratmaat}
     \left|\zeta_q(2) - \frac{a_{n}}{b_{n}}\right|
    = \mathcal{O} \left( \frac{ p^{\left(\frac{-3(\pi^2-8)}{\pi^2}+\epsilon\right)n^2}}{b_{n}}\right)
\end{equation}
for every $\epsilon > 0$, with $b_{n}=d_{2n-1}^2(p) p_{n,n-1}(p^{2n-1})$. 
Combining the result of Lemma \ref{dn-asymptotiek}
with the asymptotic behavior of $p_{n,n-1}(p^{2n-1})$ given in
(\ref{asymptotiekp}) gives 
\[
    \lim_{n \to \infty} \left| b_n \right|^{1/n^2}
    = \lim_{n \to \infty} \left| d_{2n-1}^2(p) p_{n,n-1}(p^{2n-1})\right|^{1/n^2}
    =p^{24/\pi^2}p^{6}
    =p^\frac{6(\pi^2+4)}{\pi^2}.
\]
Together with (\ref{irratmaat}) this gives
\begin{equation}
    \left|\zeta_q(2) - \frac{a_{n}}{b_{n}}\right|
    = \mathcal{O} \left( \frac{1}{ b_n^{1+ \frac{\pi^2-8}{2 (\pi^2+4)}-\epsilon}}\right)
\end{equation}
for every $\epsilon >0$, which implies the following upper bound
for the measure of irrationality
\[
    r(\zeta_q(2)) \leq 1+ \frac{2(\pi^2+4)}{\pi^2-8}= \frac{3 \pi^2}{\pi^2-8} \approx 15.8369.
\]

\section{Linear independence of $1,\zeta_q(1), \zeta_q(2)$ over $\mathbb{Q}$}
In this section we prove Theorem \ref{thm:Linind} by using Lemma
\ref{lem:Linind}. To apply this lemma we need rational
approximations for $\zeta_q(1)$ and $\zeta_q(2)$ with common
denominator. From the equations (\ref{beta})-(\ref{alpha}) and
(\ref{bnm})-(\ref{anm}), it follows that we can take
$p_n^*=b_n=d_{2n-1}(p)\beta_n$, $q_n^*=d_{2n-1}(p) \alpha_n$ and
$r_n^*=a_n$, i.e.
\begin{eqnarray}\label{pn*qn*rn*}
  \hspace*{-10pt}  p_n^* &=& d_{2n-1}^2(p)p_{n,n-1}(p^{2n-1})\\
  \hspace*{-10pt} q_n^* &=& d_{2n-1}^2(p)\left[q_{n,n-1}(p^{2n-1})+
    p_{n,n-1}(p^{2n-1})\sum_{k=1}^{2n-2}\frac{1}{p^k-1}\right]\\
  \hspace*{-10pt} r_n^* &=&
    d_{2n-1}^2(p)\left[r_{n,n-1}(p^{2n-1})+
    p_{n,n-1}(p^{2n-1})\sum_{k=1}^{2n-2}\frac{k}{p^k-1}+(2n-1)q_{n,n-1}(p^{2n-1})\right].
\end{eqnarray}
Then Theorems \ref{thm:qzeta1} and
\ref{thm:qzeta2} imply that $p_n^*$, $q_n^*$ and $r_n^*$ are
integers. The following verification of the 3 conditions of Lemma
\ref{lem:Linind} will establish the linear independence of $1$,
$\zeta_q(1)$ and $\zeta_q(2)$ over $\mathbb{Q}$. 

\subsection{Condition 1}
We will now show that $a p_n^*+b q_n^*+c r_n^* \neq 0$ for all $n
\in \mathbb{N}$ for which $2n-1$, is prime and
$2n-1>c$. It is sufficient to prove that
\[   a p_n^*+b q_n^*+c r_n^* \not \equiv 0 \mod    \Phi_{2n-1}(p) \]
or
\[     \Phi_{2n-1}(p) \nmid a p_n^*+b q_n^*+c r_n^* \]
for these values of $n$. We prove this in three steps.

\textbf{Step 1}: In the first step we prove that
\begin{equation}\label{step1}
    a p_n^*+b q_n^*+c r_n^*
    \equiv - c \, \frac{d_{2n-1}^2(p)}{(p^{2n-1}-1)^2}\mod \Phi_{2n-1}(p).
\end{equation}
First note that $d_{2n-1}(p)= \Phi_{2n-1}(p) d_{2n-2}(p)$, which implies that 
$\Phi_{2n-1}(p)$ divides the integers
$d_{2n-1}^2(p)p_{n,n-1}(p^{2n-1})$ and $d_{2n-1}^2(p)q_{n,n-1}(p^{2n-1})$
in $\mathbb{Z}$. Using this we get
\begin{eqnarray*}
  a p_n^*+b q_n^*+c r_n^*  & \equiv &  c \, d_{2n-1}^2(p)r_{n,n-1}(p^{2n-1})  \mod \Phi_{2n-1}(p)\\
    & \equiv & c \, d_{2n-1}(p)d_{2n-2}(p)\Phi_{2n-1}(p) r_{n,n-1}(p^{2n-1}) \mod  \Phi_{2n-1}(p).
\end{eqnarray*}
Expression (\ref{rnm ge-evalueerd}) implies that all the terms in
the sum $c \, d_{2n-1}(p)d_{2n-2}(p)r_{n,n-1}(p^{2n-1})$ are integers
except for the one with $r=i=2n-1$. So we see
\[
     a p_n^*+b q_n^*+c r_n^*
     \equiv - c \, \qbinom{3n-2}{n-1}_p
    \qbinom{3n-1}{n}_p
      \frac{d_{2n-1}^2(p)}{(p^{2n-1}-1)^2}\mod \Phi_{2n-1}(p).
\]
Now we only have to eliminate the binomial numbers. This can be
done by the following lemma. By applying it twice, with $(n,m)$ replaced by $(n,n-1)$ and $(n-1,n)$, we obtain
equation (\ref{step1}).
\begin{lemma}
The following congruence for polynomials in $\mathbb{Z}[x]$ holds
\begin{equation}\label{Jonathan}
       \qbinom{n+2m}{m}_x \equiv 1 \mod \Phi_{n+m}(x)
\end{equation}
for all $n, m \in \mathbb{N}$.
\end{lemma}
\begin{proof}
We prove this result by induction on $m$. Obviously relation
(\ref{Jonathan}) is satisfied for all $n \in \mathbb{N}$ when
$m=0$. Suppose that
\begin{equation}\label{inductiehypothese}
    \qbinom{n+2m-2}{m-1}_x \equiv 1 \mod \Phi_{n+m-1}(x)
\end{equation}
for all $n \in \mathbb{N}$, then we can prove 
(\ref{Jonathan}) as follows. 
The $q$-version of Pascal's triangle identity (\ref{Pascal})
gives that
\[
    \qbinom{n+2m}{m}_x=\qbinom{n+2m-1}{m}_x+x^{n+m}\qbinom{n+2m-1}{m-1}_x.
\]
Note that (\ref{x^n-1}) implies that $x^{n+m}
\equiv 1 \mod  \Phi_{n+m}(x)$. We can also write, using
(\ref{x^n-1}),
\[
    \qbinom{n+2m-1}{m}_x=
    \frac{(x;x)_{n+2m-1}}{(x;x)_{m}(x;x)_{n+m-1}}
    =\frac{\displaystyle\prod_{\nu=n+m}^{n+2m-1}(1-x^\nu)}{\displaystyle \prod_{\nu=1}^{m}(1-x^\nu)}
    = \frac{\displaystyle \prod_{\nu=n+m}^{n+2m-1}\prod_{d|\nu} \Phi_{d}(x)}
    {\displaystyle \prod_{\nu=1}^{m}\prod_{d|\nu} \Phi_{d}(x)}.
\]
Since $\qbinom{n+2m-1}{m}_x$ is a polynomial in $x$ with integer
coefficients, the cyclotomic polynomials $\Phi_{d}(x)$ are
irreducible over $\mathbb{Q}$, and $\Phi_{n+m}(x)$ is a
factor of the numerator and not of the denominator, it follows
that $\Phi_{n+m}(x)$ divides $\qbinom{n+2m-1}{m}_x$. Using
this, it turns out that
\[
    \qbinom{n+2m}{m}_x \equiv \qbinom{n+2m-1}{m-1}_x \mod \Phi_{n+m}(x).
\]
Now we can apply the induction hypothesis
(\ref{inductiehypothese}) with $n+1$ instead of $n$. This proves
the lemma.
\end{proof}

\textbf{Step 2}: In the second step we prove that
\begin{equation}\label{step2}
    \gcd\left(c,\Phi_{n+m}(p)\right)=1.
\end{equation}
For this we need the following result of Legendre (see, e.g. \cite{Gallot})
\begin{lemma}\label{Legendre}
For all positive integers $p, s$ and every cylotomic polynomial
$\Phi_{n}$ we have that if $s \mid \Phi_{n}(p)$ then $s=n\ell
+1$ for some $\ell\in \mathbb{N}_0$, or $s \mid n$.
\end{lemma}

Suppose $\gcd\left(c,\Phi_{2n-1}(p)\right)=s>1$. Then $s \mid
\Phi_{2n-1}(p)$ so from Lemma \ref{Legendre} it follows that
$s=\ell(2n-1) +1$ for some $\ell\in \mathbb{N}_0$ or $s \mid 2n-1$.
Since $s>1$ and $2n-1$ is a prime number, we have in both cases that
$s \geq 2n-1$. But $s \mid c$ and we only consider values of $n$
for which $2n-1 > c$. This gives a contradiction and proves
equation (\ref{step2}).

\textbf{Step 3}: In the last step we prove that
\begin{equation}\label{step3}
    \Phi_{2n-1}(p) \nmid \frac{d_{2n-1}^2(p)}{(p^{2n-1}-1)^2}.
\end{equation}

We will do this by contraposition: suppose that there exist an
integer $A \neq 0$ such that
\[
    \frac{d_{2n-1}^2(p)}{(p^{2n-1}-1)^2}= A \,  \Phi_{2n-1}(p)
\]
or, when we use (\ref{dn})--(\ref{x^n-1}),
\[
    \prod^{2n-1}_{\substack{k=1\\k \nmid 2n-1}}\Phi_{k}^2(p)= A \,
    \Phi_{2n-1}(p).
\]
Suppose $s$ is a prime number so that $s \mid \Phi_{2n-1}(p)$, then it
follows from the previous equation that there exists an integer
$k$, with $2\leq k \leq 2n-2$, such that $s \mid \Phi_{k}(p)$. We
first prove that this implies $\gcd(\Phi_{2n-1}(x), \Phi_{k}(x))=1$ in
$(\mathbb{Z}/s\mathbb{Z})[x]$. For this we need the following
lemma \cite[Section 4.9]{Stillwell}:
\begin{lemma}\label{lem:factor}
If $s$ is a prime not dividing $n$, then $x^n-1$ has no repeated
factors in $(\mathbb{Z}/s\mathbb{Z})[x]$.
\end{lemma}
Since $s \mid \Phi_{2n-1}(p)$ we can first apply Lemma
\ref{Legendre} and distinguish between two cases.
\begin{description}
\item[\ case 1:] $s=\ell(2n-1)+1$ for some $\ell\in \mathbb{N}_0$.  In this
case it is obvious that
 $s \nmid 2n-1$, and also
$s \nmid  k$ because $k \leq 2n-2<s$. So we can apply Lemma
\ref{lem:factor} with parameter $k(2n-1)$ instead of $n$. We get
that $x^{k(2n-1)}-1$ has no repeated factors in
$(\mathbb{Z}/s\mathbb{Z})[x]$. This means, using (\ref{x^n-1}),
that
\[
    \Phi_{k}(x)\cdot \Phi_{2n-1}(x)\cdot \prod_{\substack{d \mid k(2n-1)\\d \neq k, \, d \neq 2n-1}} \Phi_{d}(x)
\]
has no repeated factors, so $\gcd(\Phi_{2n-1}(x), \Phi_{k}(x))=1$
in $(\mathbb{Z}/s\mathbb{Z})[x]$.\\
\item[case 2:] $s \mid 2n-1$. Since $s$ and $2n-1$ are both prime, we see that
$s=2n-1$ in this case. Since $s$ is prime, we also know that
$\Phi_{s}(x)=(x^s-1)/(x-1).$ In $(\mathbb{Z}/s\mathbb{Z})[x]$
this implies
\[
    \Phi_{s}(x)=\frac{x^s-1}{x-1}=\frac{(x-1)^s}{x-1}=(x-1)^{s-1}.
\]
So $\gcd(\Phi_{k}(x), \Phi_{2n-1}(x))=\gcd(\Phi_{k}(x),
(x-1)^{s-1})=(x-1)^u$, for some $u$ with $0\leq u \leq s-1$ in
$(\mathbb{Z}/s\mathbb{Z})[x]$.
We know that $s$ is a prime not dividing $k$, so applying Lemma
\ref{lem:factor} gives that $x^k-1$ has no repeated factors in
$(\mathbb{Z}/s\mathbb{Z})[x]$. Now we can write $x^k-1$ as
\[
    x^k-1=\Phi_{1}(x)\cdot \Phi_{k}(x)\cdot \prod_{\substack{d \mid k\\1<d<k}}
    \Phi_{d}(x).
\]
and since $\Phi_{1}(x)=x-1$, it follows that $x-1$ can not be a
factor of $\Phi_{k}(x)$, so we can conclude that
$\gcd(\Phi_{2n-1}(x), \Phi_{k}(x))=1$ in
$(\mathbb{Z}/s\mathbb{Z})[x]$.
\end{description}

So in both cases we have $\gcd(\Phi_{2n-1}(x), \Phi_{k}(x))=1$ in
$(\mathbb{Z}/s\mathbb{Z})[x]$. Integer division gives us two
polynomials $u(x)$ and $v(x)$ such that $u(x)\Phi_{2n-1}(x)+
v(x)\Phi_{k}(x)=1$ in $(\mathbb{Z}/s\mathbb{Z})[x]$. Evaluating
this expression at $x=p$ gives in $\mathbb{Z}$ the
equality
\[
    u(p)\Phi_{2n-1}(p)+ v(p)\Phi_{k}(p)=1+ s  w,
\]
where $w$ is an integer. But from $s \mid \Phi_{n+m}(p)$ and $s
\mid \Phi_{k}(p)$ we see that $s \mid 1$, which gives a
contradiction and proves (\ref{step3}).
\medskip

From these three steps it follows that $\Phi_{2n-1}(p) \nmid a
p_n^*+b q_n^*+c r_n^*$ for infinity many values of $n$, so the
first condition of Lemma \ref{lem:Linind} is satisfied.

{\bf Remark:} For the case $c=0$, $b \neq 0$ we can work analogously
by considering
\[    a p_n^*+b q_n^*+c r_n^*  \mod d_{2n-1}(p)\Phi_{2n-1}(p).
\]
 The case $c=b=0$ is obvious.

\subsection{Condition 2 and 3}
From the definition of $p_n^*$ and $q_n^*$, it follows that
$|p_n^* \zeta_q(1) - q_n^*|= d_{2n-1}(p)|\beta_n \zeta_q(1) -
\alpha_n|$. By using Theorem \ref{thm:qzeta1} together with Lemma
\ref{dn-asymptotiek}, it turns out that
\[
    |p_n^* \zeta_q(1) - q_n^*|^{1/n^2}=\left(d_{2n-1}(p)\right)^{1/n^2}|\beta_n
    \zeta_q(1) - \alpha_n|^{1/n^2}=p^{-\frac{3(\pi^2-8)}{4\pi^2}}<1,
\]
so that $|p_n^* \zeta_q(1) - q_n^*|\to 0$ when $n \to \infty$.

It also follows from the definition of $p_n^*$ and $r_n^*$ that
$|p_n^* \zeta_q(2) - r_n^*|  = |b_n \zeta_q(2) - a_n|$ and from
Theorem \ref{thm:qzeta2} we know already that this expression
tends to zero when $n$ tends to infinity.


\begin{thebibliography}{99}
\bibitem{AAR} G. E. Andrews, R. Askey, R. Roy,
    \textit{Special Functions}, Encyclopedia of Mathematics and its
    Applications \textbf{71}, Cambridge University Press, 1999.
\bibitem{Ap} R. Ap\'ery,
    \textit{Irrationalit\'e de $\zeta(2)$ et $\zeta(3)$},
    Ast\'erisque \textbf{61} (1979), 11--13.
\bibitem{BallRivoal} K. Ball, T. Rivoal,
    \textit{Irrationalit\'e d'une infinit\'e de valeurs de la
    fonction z\`eta aux entiers impairs},
    Invent. Math. \textbf{146} (2001), no.~1, 193--207.
\bibitem{Borwein-Borwein}
    J.M. Borwein and P.B. Borwein,
    \textit{Pi and the AGM--A Study in Analytic Number Theory
    and Computational Complexity}, Wiley, New York, 1987.
\bibitem{Driver-Stahl}
    K. Driver and H. Stahl,
    \textit{Simultaneous rational approximants to Nikishin-systems I},
    Acta Sci. Math. \textbf{60} (1995), 245--263.
\bibitem{Fischler}
    S. Fischler,
    \textit{Irrationalit\'e de valeurs de z\^eta (d'apr\`es Ap\'ery, Rivoal,\ldots)},
    Ast\'erisque \textbf{294} (2004), 27--62.
\bibitem{Gallot} Y. Gallot,
    \textit{Cyclotomic polynomials and prime numbers}, January 2001, \\
    \texttt{http://perso.wanadoo.fr/yves.gallot/papers/cyclotomic.html}
\bibitem{GR} G. Gasper, M. Rahman,
    \textit{Basic Hypergeometric Series},
    Encyclopedia of Mathematics and its Applications \textbf{35},
    Cambridge University Press, 1990.
\bibitem{Hardy-Wright}
    G.H. Hardy, E.M. Wright,
    \textit{An Introduction to the Theory of Numbers},
    Oxford University Press, 1938 (5th edition, 1979).
\bibitem{Hata} M. Hata,
    \textit{Rational approximations to $\pi$ and some other
    numbers},
    Acta Arith. \textbf{63} (1993), no.~4, 335--349.
\bibitem{KRZ} C. Krattenthaler, T. Rivoal, W. Zudilin,
    \textit{S\'eries hyperg\'eom\'etriques basiques, $q$-analogues des valeurs de
    la fonction z\^eta et s\'eries d'Eisenstein},
    J. Inst. Math. Jussieu \textbf{5} (2006), 53--79.
\bibitem{Nesterenko} Yu. V. Nesterenko,
    \textit{Modular functions and transcendence questions},
    Mat. Sbornik \textbf{187} (1996), no.~9, 65--96 (in Russian);
    Sbornik Math. \textbf{187} (1996), no.~9, 1319--1348.
\bibitem{NikiSor} E. M. Nikishin, V. N. Sorokin,
    \textit{Rational Approximations and Orthogonality},
    Translations of Mathematical Monographs \textbf{92},
    Amer. Math. Soc., Providence RI, 1991.
\bibitem{Walter-Kelly} K. Postelmans, W. Van Assche,
    \textit{Multiple little $q$-Jacobi polynomials},
    J. Comput. Appl. Math. \textbf{178} (2005), 361--375.
\bibitem{Ribenboim} P. Ribenboim,
    \textit{The New Book on Prime Number Records},
    Springer-Verlag, New York, 3rd edition 1995.
\bibitem{Rivoal} T. Rivoal,
    \textit{Irrationalit\'e d'au moins un des neuf nombres
    $\zeta(5),\zeta(7),\ldots,\zeta(21)$},
    Acta Arith. \textbf{103} (2002), no.~2, 157--167.
\bibitem{Stillwell} J. Stillwell,
    \textit{Elements of Algebra},
    Springer-Verlag, New York, 1996.
\bibitem{VanAssche}
    W. Van Assche,
    \textit{Little $q$-Legendre polynomials and irrationality of certain
    Lambert series},
    Ramanujan J. \textbf{5} (2001), 295--310.
\bibitem{VanAssche2} W. Van Assche,
    \textit{Multiple orthogonal polynomials, irrationality
    and transcendence},
    in `Continued Fractions: from analytic number
    theory to constructive approximation' (B. C. Berndt and F. Gesztesy, eds.),
    Contemporary Mathematics \textbf{236}, Amer. Math. Soc., Providence RI, 1999,
    pp.~325--342.
\bibitem{Zudilin3}
    V. V. Zudilin,
    \textit{One of the numbers
    $\zeta(5),\zeta(7),\zeta(9),\zeta(11)$ is irrational},
    Uspekhi Mat. Nauk \textbf{56} (2001), no.~4, 149--150 (in
    Russian);
    Russian Math. Surveys \textbf{56} (2001), no.~4,
    774--776.
\bibitem{Zudilin}
    V. V. Zudilin,
    \textit{Irrationality of  values of the Riemann zeta
    function},
    Izv. Ross. Akad. Nauk Ser. Mat. \textbf{66} (2002), no.~3,
    49--102 (in Russian);
    Izv. Math. \textbf{66} (2002), no.~3 (2002),
    489--542.
\bibitem{Zudilin1}
    V. V. Zudilin,
    \textit{On the irrationality measure of the $q$-analogue of $\zeta(2)$},
    Mat. Sbornik \textbf{193} (2002), no.~8, 49--70 (in Russian);
    Sbornik Math. \textbf{193} (2002), no.~7-8, 1151--1172.
\bibitem{Zudilin2}
    V. V. Zudilin,
    \textit{Diophantine problems for $q$-zeta values},
    Mat. Zametki \textbf{72} (2002), no.~6, 936--940 (in Russian);
    Math. Notes \textbf{72} (2002), no.~5--6, 858--862.

\end{thebibliography}
\end{document}